\documentclass[12pt,envcountsect]{svjour3} 
\smartqed 
\usepackage{graphicx}
\usepackage{caption}
\usepackage{booktabs,multirow}
\usepackage{amssymb,amsmath}
\usepackage[all]{xy}
\usepackage{float}
\usepackage{subfig}
\usepackage{graphicx}

\usepackage{framed} 
\usepackage{lastpage}

\usepackage[algoruled]{algorithm2e}

\usepackage[pagebackref=false,colorlinks,linkcolor=blue,citecolor=magenta]{hyperref}
\usepackage[nottoc]{tocbibind}
\journalname{JOTA}

\newtheorem{assumption}{Assumption}[section]

\begin{document}

\title{A Gradient Sampling method based on Ideal direction for solving nonsmooth nonconvex optimization problems}
\subtitle{ convergence analysis and numerical experiments}

\author{M. Maleknia   \and  M. Shamsi }

\institute{ Morteza Maleknia \at
				Amirkabir University of Technology\\ 
				Tehran, Iran\\
    			m.maleknia@aut.ac.ir  
           \and
           		M. Shamsi,  Corresponding author  \at
				Amirkabir University of Technology\\ 
				Tehran, Iran\\
               m\_shamsi@aut.ac.ir
}

\date{Received: date / Accepted: date}

\maketitle

\begin{abstract}
In this paper, a modification to the Gradient Sampling (GS) method for minimizing  nonsmooth nonconvex functions is presented. 	 
One drawback in GS method is the need of solving a Quadratic optimization Problem (QP) at each iteration, which is time-consuming especially for large scale objectives.
To resolve this difficulty, we propose a new descent direction, namely Ideal direction, for which there is no need to consider any quadratic or linear optimization subproblem. It is shown that, this direction satisfies Armijo step size condition and can be used to make a substantial reduction in the objective function. Furthermore, we prove that using Ideal directions preserves the global convergence of the GS method. Moreover, under some moderate assumptions, we present an upper bound for the number of serious iterations. Using this upper bound, we develop a different strategy to study the convergence of the method. We also demonstrate the efficiency of the proposed method using small, medium and large scale problems in our numerical experiments.
\end{abstract}


\newcommand{\I}{ {\text{\rm\tiny{I}}} }
\newcommand{\s}{ {\text{\rm\tiny{s}}} }

\section{Introduction}
In this paper, we consider the following unconstrained minimization problem
\begin{equation}\label{Main-Problem}
\min f({\mathbf x}) \quad \text{s.t.} \quad {\mathbf x}\in\mathbb{R}^n 
\end{equation}
where $f:\mathbb{R}^n\to\mathbb{R}$ is locally Lipschitz and continuously differentiable on an open set $D$ with full measure in $\mathbb{R}^n$. 
These types of problems arise in many applications, such as optimal control, image processing, data analysis, economics, chemistry and biology (\cite{Makela_book,image,Poly-approx,bio-ch2,data-mining}). Therefore, it is worthwhile to develop efficient algorithms for solving such problems. During the last three decades, a lot of effort has gone into nonsmooth optimization. In particular, Clarke developed the concept of subgradient in \cite{Clarke1990} and soon after, it became the cornerstone of designing and developing many algorithms in nonsmooth optimization (see \cite{Makela_book,kiwielbook,shorbook,lemarchal-book1,lemarchal-book2,v-metric1,v-metric2}).

   When it comes to designing numerical methods for minimizing nonsmooth functions, there are serious challenges we need to deal with. The main difficulty is that, nonsmooth functions are generally not differentiable at stationary points. Finding a descent direction is not an easy task as well. In smooth case, any vector has an obtuse angle with the gradient is a descent direction. In particular, the vector $-\nabla f({\mathbf x})$ defines the steepest descent direction. In contrast, in the nonsmooth case a vector that is opposite to an arbitrary subgradient need not be a direction of descent and hence designing descent algorithms is much more complicated. These difficulties make most of classic methods in smooth optimization unsuitable for solving nonsmooth optimization problems. For instance, it is well known that when  the ordinary steepest descent method is applied to a nonsmooth function, it generally fails to make a substantial reduction in objective function as it approaches a nonsmooth region.

There are various methods for locating  minimizers of nonsmooth functions. The subgradient method, originally proposed by N. Shor \cite{shorbook}, is one of the simplest methods for minimizing nonsmooth functions. Due to its simple structure, it is a popular  method, although it suffers from some serious limitations such as slow convergence, lack of practical termination criterion and lack of descent. However, it is well known that some of its modifications are able to overcome these difficulties (see \cite{bagirov2012,m-amiri,Nesterov,bagirov2010}). As another class of nonsmooth methods, we can refer to Bundle methods as one of the most efficient methods in nonsmooth optimization \cite{lemarchal-book1,lemarchal-book2,kiwielbook}. The key idea of bundle methods is to keep memory of computed subgradients at previous iterations to construct a piecewise linear model for objective function. These methods require a great deal of storage, therefore some of its modifications presented in \cite{Haarala,Luksan} are more efficient for solving large scale nonsmooth problems. In class of subgradient and bundle methods, at each iteration, the user needs to supply at least one subgradient. However, in many cases computing only one subgradient is not an easy task. In such situations, Gradient Sampling (GS) methods are more efficient, because they obtain a search direction without explicit computation of subgradients.

The gradient sampling algorithm, originally developed by Burke, Lewis and Overton \cite{Burke2005}, is a descent method for solving problem \eqref{Main-Problem}. The method is robust and can be applied to a broad range of nonsmooth functions. A comprehensive discussion of the GS algorithm along with its last modifications can be found in the recent paper \cite{GS-full}. The convergence analysis and theoretical results of the GS algorithm were first developed in \cite{Burke2005}. Soon after, these results were strengthen by the work of Kiwiel in \cite{Kiwiel2007}. In particular, under subtle modifications, Kiwiel derived a lower bound for step sizes and suggested a limited Armijo line search in which the number of backtracking steps can be managed through our choice of initial step size. Another version of the GS method in which the Clarke $\varepsilon$-subdifferential is approximated by sampling estimates of mollifier gradients presented by Kiwiel in \cite{Kiwiel2010}. In the work of Curtis and Que \cite{Curtis2013}, two novel strategies were introduced to improve the efficiency of the GS method. The first strategy is to combine the idea of sampling gradients and LBFGS update (see \cite{nocedal}) to approximate Hessian matrix. The second one uses the idea of sampling gradients to provide a model that overestimates the objective function. In both techniques, the dual  of the quadratic subproblem is considered for warm-starting the QP solver. In the work of Curtis and Overton \cite{Curtis2012}, the ideas of gradient sampling and Sequential Quadratic Programming (SQP) are combined for solving nonsmooth constrained optimization problems. In addition, the local convergence rate of the GS algorithm for the class of finite-max functions is studied in \cite{C-Rate}. These continued developments clearly indicate that, the GS technique has been a reach area of research.

However, the GS method suffers from two practical limitations. First, in order to obtain an approximation of  $\varepsilon$-steepest descent direction, the GS method requires to compute gradient information at $m\geq n+1$ randomly generated points during each iteration. It is noted that, for large scale objectives, this may not be tractable at a moderate cost. To alleviate this difficulty, In \cite{Curtis2013},  Curtis and Que proposed an adaptive strategy in which the convergence of the method is guaranteed only through   $\mathcal{O}(1)$ gradient evaluations at each iteration. The second limitation is to solve the corresponding QP. As the number of variables increases, the size of this QP increases significantly which makes the method unsuitable for large scale objectives. Furthermore, when we are far away from a nonsmooth region, since the information  collected by sampling gradients are very close together, solving the quadratic subproblem is not reasonable. Nevertheless, since checking differentiability of $f$ is not an easy task, the GS method and its modifications do not care about this fact. In this work, we address the second difficulty by introducing a new search direction.

In this paper, based on gradient sampling technique, we propose a new search direction, namely Ideal direction, for which there is no need to consider any kind of quadratic or linear subproblems. Furthermore, we show that this direction satisfies Armijo step size condition and provides a necessary optimality condition in the sense that it vanishes at optimality. The original GS method is modified using Ideal directions. By means of this new search direction, not only is there no need  to solve the quadratic subproblem in smooth regions, but also we can reduce the number of quadratic subproblems once the method starts tracking a nonsmooth curve which leads to a stationary point. Moreover, we provide a comprehensive discussion for the convergence analysis of the method. We follow closely the work of Kiwiel in \cite{Kiwiel2007} to analyze the convergence of the proposed method. Of course, there are some differences due to using Ideal directions. In addition, thanks to the limited Armijo line search proposed in \cite{Kiwiel2007}, we present an upper bound for the number of serious iterations generated in our method. Using this upper bound, a different strategy to study the global convergence behavior of the proposed method is developed assuming that the objective function is bounded below.

This paper is organized as follows. Section \ref{section2} provides some mathematical preliminaries used in this paper. A brief review over the GS method is presented in Sect. \ref{Section3}. In Sect. \ref{Section4}, we introduce Ideal direction and its main properties are examined. The proposed method and its convergence analysis are given in Sect. \ref{Section5}. Numerical results are reported in Sect. \ref{Section6} and Sect. \ref{Conclusion} concludes the paper. 

\section{Preliminaries}\label{section2}

We use the following notations in this paper. As usual, $\mathbb{R}^n$ is the $n$-dimensional Euclidean space and its inner product is denoted by $\langle {\mathbf x} ,\mathbf y\rangle := \sum_{k=1}^{n} x_i y_i$ that induces the associated Euclidean norm $\lVert {\mathbf x} \rVert$:=$\langle {\mathbf x} , {\mathbf x}\rangle^{1/2} $. $B({\mathbf x},\varepsilon):=\{\mathbf y\in \mathbb{R}^n : \lVert \mathbf y-{\mathbf x} \rVert \leq \varepsilon \}$ is the closed ball centered at ${\mathbf x}$ with radius $\varepsilon$. Furthermore $B_{\varepsilon}:=B(\mathbf 0,\varepsilon)$. \\

A function $f:\mathbb{R}^n\to\mathbb{R}$ is called locally Lipschitz \cite{Bagirov2014}, if  for every ${\mathbf x}\in\mathbb{R}^n$ there exist positive constants $K_{\mathbf x}$ and $\varepsilon_{\mathbf x}$ such that   
$$ \vert f(\mathbf z)-f(\mathbf y)\vert \leq K_{\mathbf x} \lVert \mathbf z-\mathbf y \rVert \, , \quad \text{for \,all}\,\, \mathbf z ,\mathbf y \in B({\mathbf x},\varepsilon_{\mathbf x}). $$   
Let
$$\Omega_f:=\{{\mathbf x}\in\mathbb{R}^n \, : \, f\,\, \text{is not differentiable at } {\mathbf x}\}, $$
be the subset of $\mathbb{R}^n$ where the function $f$ is not differentiable. By Rademacher's theorem \cite{Evans2015}, every locally Lipschitz function is differentiable almost everywhere. Therefore, the Clarke subdifferential of a locally Lipschitz function  $f$ at a point ${\mathbf x}\in\mathbb{R}^n$ can be given by \cite{Clarke1990}
$$\partial f({\mathbf x}):=\texttt{co}\{{\boldsymbol \xi}\in \mathbb{R}^n : \exists \, \{{\mathbf x}_k\}\subset \mathbb{R}^n\setminus\Omega_f  \,\,\,\, \text{s.t.} \,\,\, \, {\mathbf x}_k\to {\mathbf x}\,\,\, {\rm and} \,\,\,\, \nabla f({\mathbf x}_k)\to {\boldsymbol \xi} \},$$
in which $\texttt{co}$ stands for the convex hull. It is shown in \cite{Bagirov2014} that the set valued map $\partial f:\mathbb{R}^n\rightrightarrows \mathbb{R}^n$ is outer semicontinuous and the set $\partial f({\mathbf x})$ is a convex compact subset of $\mathbb{R}^n.$ The Clarke $\varepsilon$-subdifferential, which is the generalization of ordinary subdifferential, is defined by \cite{Bagirov2014}
$$\partial_{\varepsilon}f({\mathbf x}):=\texttt{cl\,co}\, \partial f({\mathbf x}+B_{\varepsilon}) .$$
Clearly, the set valued map $\partial_{\varepsilon} f:\mathbb{R}^n\rightrightarrows \mathbb{R}$ has a closed graph and hence it is outer semicontinuous \cite{Rockafellar2004}. \\

Suppose that $f:\mathbb{R}^n\to\mathbb{R}$ is locally Lipschitz and continuously differentiable on an open set $D$ with full measure in $\mathbb{R}^n$. The following representation of the Clarke subdifferential is the key idea of approximating subdifferential set by sampling gradients \cite{Clarke1990}
$$\partial f({\mathbf x})=\bigcap _{\varepsilon>0} G_{\varepsilon}({\mathbf x}) ,$$ 
such that
$$G_{\varepsilon}({\mathbf x}):= \texttt{cl\,co} \{\nabla f(({\mathbf x}+B_{\varepsilon})\cap D)\},$$
where $\texttt{cl}$ denotes the closure of a set. Let ${\mathbf x}\in\mathbb{R}^n$ be a differentiable point, $\varepsilon >0$, $m\in\mathbb{N}$ and $\mathbf u_1,...,\mathbf u_m$ be sampled uniformly and independently from $B_1$. If $\mathbf s_0:={\mathbf x}$ and $\mathbf s_i:={\mathbf x}+\varepsilon \mathbf u_i$ is a differentiable point for $i=1,\ldots,m$, then the gradient bundle $G^m_\varepsilon(\mathbf x)$ is defined as \cite{Burke2005,Kiwiel2007}  
\begin{equation}\label{Subdifferential approximation}
G^m_\varepsilon({\mathbf x}):= \texttt{co} \{\nabla f(\mathbf s_0),\nabla f(\mathbf s_1), \ldots, \nabla f(\mathbf s_m) \},
\end{equation}
in which $\epsilon>0$ and $m\in\mathbb{N}$ are called sampling radius and sample size respectively.
For a locally Lipschitz function $f$, it is easy to see that 
$$G^m_\varepsilon({\mathbf x}) \subset G_\varepsilon({\mathbf x}) \subset \partial_{\varepsilon} f({\mathbf x}) ,$$
and for $0<\varepsilon_1<\varepsilon_2$, we have $\partial_{\varepsilon_1}f({\mathbf x}) \subset G_{\varepsilon_2}({\mathbf x})$. Furthermore, it is proved in \cite{GS-full}   that for every ${\mathbf x}\in\mathbb{R}^n$ and $\varepsilon>0$ 
$$\partial f({\mathbf x})\subset G_\varepsilon({\mathbf x}). $$
The following definition provides a useful point of view to study the concept of gradient sampling.
\begin{definition}
	For a locally Lipschitz function $f:\mathbb{R}^n\to\mathbb{R}$ and optimality tolerance $\nu>0$, a point $\mathbf x\in\mathbb{R}^n$ is called a $(\nu,G^m_\epsilon(\mathbf x))$-stationary point if
	$$\min \{\lVert \mathbf g \rVert : \mathbf g\in G^m_\epsilon(\mathbf x) \} \leq \nu .$$
\end{definition}

We recall that, if ${\mathbf x}\in\mathbb{R}^n$ is a local minimum of a locally Lipschitz function $f$, then it is necessary that \cite{Bagirov2014}

$$\mathbf 0\in \partial f({\mathbf x}).$$
The point ${\mathbf x}\in\mathbb{R}^n$ satisfying the above condition is called Clarke stationary point. Furthermore, a point ${\mathbf x}\in\mathbb{R}^n$ is called Clarke $\varepsilon$-stationary point if
$$\mathbf 0\in \partial_\varepsilon f({\mathbf x}) .$$

\section{{ Background on GS method}}\label{Section3}
The Gradient Sampling (GS) method \cite{Burke2005} is a descent method  for solving problem \eqref{Main-Problem}. At each iteration of the method, the steepest descent direction is approximated and utilized to make a substantial reduction in $f$. In this method, the concepts of steepest descent direction and choosing step size are crucial. In this respect, these concepts are briefly reviewed .

\subsection{ The steepest descent direction}\label{The steepest descent direction}

We recall that, a direction ${\mathbf d}\in\mathbb{R}^n$ is called a direction of descent for $f:\mathbb{R}^n\to\mathbb{R}$ at ${\mathbf x}\in\mathbb{R}^n$, if there is $z_0>0$ such that
\begin{equation}\label{Descent direction}
f({\mathbf x}+t\mathbf d)-f({\mathbf x})<0 , \qquad \text{for \,all} \,\, t\in (0,z_0]. 
\end{equation}
In the descent methods, choosing a step size $t$ satisfying the above inequality  does not guarantee the convergence of the method. In this regard,
some sufficient decrease conditions are introduced \cite{kiwielbook}. 
Among them, an Armijo step size condition, which is commonly used in GS methods, is considered to state the following definition. 

\begin{definition}\label{ADD_Definition}
	Let $\mathbf g\neq \mathbf 0$, we say that $\mathbf g$ induces an Armijo Descent Direction (ADD) $\mathbf d:=-\lVert \mathbf g \rVert^{-1} \mathbf g $ for $f$ at point $\mathbf x$, if for  each $c\in(0,1)$ there is $z_0>0$ such that
	$$ f({\mathbf x}+t \mathbf d  )-f({\mathbf x})<- c t\lVert \mathbf g \rVert, \qquad {\rm{for \,\,all}} \,\, t\in (0,z_0]. $$	
	The above  condition is called \emph{sufficient decrease condition}. 
	
\end{definition}
If $f$ is smooth at a point ${\mathbf x}\in\mathbb{R}^n$, then every descent direction at this point is an ADD too\cite{nocedal}. However, this is not the case when $f$ is not smooth at ${\mathbf x}$. Therefore, in order to design a descent method for minimizing a nonsmooth function, providing an ADD is essential to ensure the convergence of the method.

To define the concept of steepest descent direction for a locally Lipschitz function more precisely, we need to generalize the classical directional derivative. For a smooth or convex function $f:\mathbb{R}^n\to\mathbb{R}$, the classical directional derivative is defined as 
\begin{align}\label{Classical directional derivative}
f'({\mathbf x};{\mathbf d}):=\lim_{t\downarrow 0}t^{-1}[f({\mathbf x}+t \mathbf d)-f({\mathbf x})],
\end{align}
and the steepest descent direction at the point ${\mathbf x}\in \mathbb{R}^n$ is obtained by solving the problem
$$\min_{\lVert {\mathbf d} \rVert \leq 1} f'({\mathbf x};{\mathbf d}) .$$
However, the quantity \eqref{Classical directional derivative} does not always exist for locally Lipschitz functions. In this regard, for locally Lipschitz functions, the Clarke generalized directional derivative  is defined as \cite{Bagirov2014}
$$f^{\circ}({\mathbf x};{\mathbf d}):=\limsup_{\substack{\mathbf y\to {\mathbf x} \\ t\downarrow 0}} \, t^{-1} [f(\mathbf y+t \mathbf d)-f(\mathbf y)] .$$ 
Moreover, it is shown that this directional derivative can be expressed as \cite{Bagirov2014}
$$f^{\circ}({\mathbf x};{\mathbf d})=\max_{{\boldsymbol \xi} \in \partial f({\mathbf x})} \langle {\boldsymbol \xi} , {\mathbf d} \rangle ,$$ 
and naturally, the steepest descent direction of a locally Lipschitz function $f$ at the point ${\mathbf x}\in\mathbb{R}^n$ can be obtained by solving the following $\min$-$\max$ problem 
\begin{equation}\label{min_max_problem}
\min_{\lVert {\mathbf d} \rVert \leq 1}\, \max_{{\boldsymbol \xi} \in \partial f({\mathbf x})} \langle {\boldsymbol \xi} , {\mathbf d} \rangle.
\end{equation}
Using von Neumann minimax theorem \cite{convex_book}, the above problem is equivalent to the following problem

\begin{equation}\label{Steepst descent direction}
\min_{{\boldsymbol \xi}\in \partial f({\mathbf x})} \, \lVert {\boldsymbol \xi} \rVert,
\end{equation} 
in the sense that, if ${\boldsymbol \xi}^*\neq 0$ is the solution of the problem \eqref{Steepst descent direction}, then ${\mathbf d}^*=-\lVert {\boldsymbol \xi}^*\rVert^{-1}{\boldsymbol \xi}^*$ solves the problem \eqref{min_max_problem}. As a result, to obtain the steepest descent direction it is sufficient to find the member of subdifferential set with minimum norm.

In an iterative descent method for minimizing a nonsmooth function, after some iterations, the method rapidly approaches a nonsmooth curve leading to a stationary point. Therefore, it is likely that, the method generates a sequence of smooth points that are close to the nonsmooth curve. In such cases, since the Clarke subdifferential is a singleton set, the steepest descent direction is not an efficient direction, while the Clarke $\varepsilon$-subdifferential is able to collect some information of the nonsmooth curve and hence the so-called $\varepsilon$-steepest descent direction is a more suitable search direction.  

In a similar fashion, the $\varepsilon$-steepest descent direction at the point ${\mathbf x} \in \mathbb{R}^n$ can be obtained from the solution of the following problem
\begin{equation}\label{epsilon Steepst descent direction}
\min_{{\boldsymbol \xi}\in \partial_{\varepsilon} f({\mathbf x})} \, \lVert {\boldsymbol \xi} \rVert.
\end{equation}
However, since solving the problem \eqref{epsilon Steepst descent direction} requires the knowledge of the whole subdifferential on $B({\mathbf x},\varepsilon)$, computing the $\varepsilon$-steepest descent direction is cumbersome. To overcome this drawback, since $G^m_\varepsilon({\mathbf x})$ is a proper inner approximation of $ \partial_\varepsilon f({\mathbf x})$, the $\varepsilon$-steepest descent direction is approximated through replacing $\partial_\varepsilon f({\mathbf x})$ by $G^m_\varepsilon({\mathbf x})$ in the problem \eqref{epsilon Steepst descent direction} and the approximate $\varepsilon$-steepest descent direction can be obtained from the following problem
\begin{equation}\label{approximate epsilon Steepst descent direction}
\min_{{\mathbf g}\in G^m_\varepsilon({\mathbf x})} \, \lVert {\mathbf g} \rVert.
\end{equation}
Indeed, the problem \eqref{approximate epsilon Steepst descent direction} is a Quadratic Problem (QP). To state this fact, by Carath\'{e}odory's theorem \cite{Bazaraa2006}, this problem can be written in the following quadratic programing form

\begin{align}\label{CQP}
 \qquad\quad&\min_{\boldsymbol \lambda } \,\frac{1}{2} \lVert G \boldsymbol \lambda \rVert^2 \\& \text{ s.t.}  \nonumber  \quad  \sum_{i=0}^{m} \lambda_i=1,\quad    \lambda_i\geq 0 ,\quad \, i=0,\ldots,m, \nonumber
\end{align}
where $G:=[\nabla f(\mathbf s_0)\,\, \nabla f(\mathbf s_1) \ldots \nabla f(\mathbf s_m)]\in\mathbb{R}^{n\times (m+1)}$ and $\boldsymbol\lambda^t:=(\lambda_0, \lambda_1, \ldots,\lambda_m)\in\mathbb{R}^{m+1}$. We summarize the concept of approximate $\varepsilon$-steepest descent direction in the following definition.
\begin{definition}
	Suppose that $\mathbf 0\notin G^m_\varepsilon({\mathbf x})$ and $\boldsymbol \lambda^*$ is the solution of the problem \eqref{CQP}. Let ${\mathbf g}^\s_\varepsilon:=G \boldsymbol\lambda^*$. Then the direction $\mathbf d^\s_\varepsilon:=-\lVert {\mathbf g}^\s_\varepsilon\rVert^{-1}  {\mathbf g}^\s_\varepsilon$ is called (normalized) approximate $\varepsilon$-steepest descent direction. 
\end{definition} 

It is noted that, the  approximate $\varepsilon$-steepest descent direction is an ADD \cite{Burke2005}. In other words, for the direction $ {\mathbf d}^\s_\varepsilon $ and each $c\in (0,1)$, there exists $z_0>0$ such that 
\begin{equation}\label{d_eps is an ADD}
f({\mathbf x}+t  {\mathbf d}^\s_\varepsilon)-f({\mathbf x})< -tc\lVert {\mathbf g}^\s_\varepsilon\rVert , \qquad \text{for \,all}\, \, t\in(0,z_0].
\end{equation}

\subsection{Choosing stepsize and backtracking line search}\label{Choosing step size}
Loosely speaking, the GS method is a descent method in which at each iteration the approximate $\varepsilon$-steepest descent direction is considered as the search direction. Therefore, the QP \eqref{CQP} is solved per iteration and then a backtracking line search is applied to find the step size $t$ as stated in Algorithm~\ref{Algorithm 1(BLS)}.

\LinesNumbered
\IncMargin{1em}
\begin{algorithm}[H]
	\SetAlgoLined

	\SetKwInOut{Input}{inputs}
	\SetKwInOut{Output}{output}
	\SetKwInOut{Required}{requied}
	
	\Indm  
	\Input{${\mathbf x}\in\mathbb{R}^n$, the directions ${\mathbf d},\mathbf g\in \mathbb R^n$, parameters   $\gamma, c\in(0,1)$. }
	\Output{Step size $t$.}
	\Required{${\mathbf d}=-\lVert \mathbf{g}\rVert^{-1} \mathbf g$ is an ADD for $f$ at ${\mathbf x}$.}
	\Indp
	\BlankLine

    \SetKwFunction{FMain}{BALS}
\SetKwProg{Fn}{Function}{:}{}
\Fn{\FMain{$\mathbf x$, $\mathbf g$, $\mathbf d$, $\gamma$, $c$}}{
		$t:=1$\;
	\While{$f({\mathbf x}+t{\mathbf d})-f({\mathbf x})\geq-c t \lVert {\mathbf g} \rVert$}{ 
		$t:=\gamma t$\;
	
	}
	\textbf{return} $t$\; 
}
\textbf{End Function}
\caption{Backtracking Armijo Line Search (BALS) }
	\label{Algorithm 1(BLS)}
\end{algorithm}
\DecMargin{1em}

In Algorithm~\ref{Algorithm 1(BLS)}, the following termination condition 
\begin{equation}\label{sufficient decrease}
f({\mathbf x}+t \mathbf d)-f({\mathbf x})< -ct\lVert \mathbf g\rVert ,
\end{equation} 
is called sufficient decrease condition. 
It is stressed that, for an ADD this algorithm finds a step size $t > 0 $ after finitely iterations.
\subsection{GS method}
The main aim of the GS method is to produce a Clarke stationary point through a sequence of $(\nu_k, G^m_{\varepsilon_k}(\mathbf x_k))$-stationary points where $\{\nu_k\}$ and $\{\varepsilon_k\}$ are decreasing sequences that tend to zero. To provide more details, at each iteration, the set of sampled points $\{\mathbf s_{k1},\ldots, \mathbf s_{km}\}\cup \{\mathbf s_{k0}\}$ is generated uniformly and independently from $B({\mathbf x}_k,\varepsilon_k)$ and then the gradient information of the objective function is computed at these points. Next,   ${\mathbf g}^\s_{\varepsilon_k}$ is obtained by solving the QP \eqref{CQP} and then the optimality condition 
\begin{equation}\label{Opt_Tol}
\lVert {\mathbf g}^\s_{\varepsilon_k}\rVert\leq \nu_k, 
\end{equation}
is checked. If condition \eqref{Opt_Tol} does not hold, ${\mathbf x}_k$ is not a $(\nu_k, G^m_{\varepsilon_k}(\mathbf x_k))$-stationary point. In this case, the sampling radius  and the optimality tolerance remain unchanged, the approximate $\varepsilon$-steepest descent direction $\mathbf d^\s_{\varepsilon_k}=-\lVert \mathbf g^\s_{\varepsilon_k}\rVert^{-1}\mathbf g^\s_{\varepsilon_k}$ is considered as the search direction, the the step size $t_k$ is computed by the Algorithm~\ref{Algorithm 1(BLS)} using inputs $(\mathbf x_k, \mathbf g^\s_{\varepsilon_k},\mathbf d^\s_{\varepsilon_k},\gamma,c)$ and finally the current point ${\mathbf x}_k$ is updated by ${\mathbf x}_{k+1}:={\mathbf x}_k+t_k \mathbf d^\s_{\varepsilon_k}$. To ensure the convergence of the GS method, the point $\mathbf x_{k+1}$ must be a differentiable point. Consequently, if the objective function is not  differentiable at ${\mathbf x}_{k+1}$, a differentiable point, say $\mathbf{\hat{x}}$, having the following properties \cite{Kiwiel2007}
\begin{align}
&f(\mathbf{\hat{x}})-f({\mathbf x}_k)<-ct_k\lVert \mathbf g^\s_{\varepsilon_k}\rVert , \label{Inequality1} \\
& \lVert {\mathbf x}_k+t_k  \mathbf d^\s_{\varepsilon_k}-\mathbf{\hat{x}}\rVert \leq \min\{t_k,\varepsilon_k\}, \label{Inequality2}
\end{align}
is selected as a perturbation of ${\mathbf x}_{k+1}$ and ${\mathbf x}_{k+1}$ is updated by $\mathbf{\hat{x}}$. A procedure to find such a perturbation can be found in \cite{Kiwiel2007}. We shall see later the motivation of the inequalities \eqref{Inequality1} and \eqref{Inequality2} in Sect. \ref{Section5}. On the other hand, if condition \eqref{Opt_Tol} holds, ${\mathbf x}_k$ is a $(\nu_k, G^m_{\varepsilon_k}(\mathbf x_k))$-stationary point. In this case, the sampling radius and the optimality tolerance are decreased by their corresponding reduction factors, the gradient bundle $G^m_{\varepsilon_k}(\mathbf x_k)$ is resampled and the process is repeated.

The GS method is robust and can be used for solving a broad range of nonsmooth problems, but it still suffers from some limitations. Here, we highlight some of these drawbacks that are addressed in this paper.
\begin{itemize}

	\item[1.] For the large scale objectives the QP \eqref{CQP} is a time-consuming subproblem that is solved at each iteration. This fact makes the GS method unsuitable for large scale nonsmooth problems \cite{GS-full}. 
	\item[2.] When the objective function $f$ is smooth on $B({\mathbf x}_k,\varepsilon_k)$ (this is the case in early iterations of the method) the information we collect by sampling gradients are close together and there is not any significant difference between $-\lVert {\mathbf g}^\s_{\varepsilon_k} \rVert^{-1} {\mathbf g}^\s_{\varepsilon_k}$ and $-\lVert \nabla f({\mathbf x}_k) \rVert^{-1} \nabla f({\mathbf x}_k)$  as a search direction. In such cases, solving the QP \eqref{CQP} is not reasonable.
\end{itemize}
To resolve the above limitations, we introduce an alternative direction, namely Ideal direction, that can be computed fast in comparison to solving the QP \eqref{CQP}. Moreover, this direction may help us to check the optimality condition \eqref{Opt_Tol} without solving the QP \eqref{CQP}.

\section{Ideal direction}\label{Section4}
To overcome the aforementioned difficulties, we introduce Ideal directions as an alternative of ${\mathbf d}^\s_{\varepsilon}$. The name ``Ideal'' originates from the concept of multicriteria optimization where Ideal points are used as reference points in compromise programming \cite{Ehrgott2005}.

For a differentiable point ${\mathbf x}\in\mathbb{R}^n, m\in\mathbb{N}$, $\varepsilon>0$ and $i\in\{1,\ldots,n\}$  let 
 $$G^m_{\varepsilon,i}(\mathbf x):=\texttt{co}\left\{\frac {\partial f}{\partial x_i}(\mathbf s_j) \,\, : \,\, j=0,\ldots,m \right\}, $$
 and define the vector $ \mathbf g^\I_\varepsilon=(g^\I_{\varepsilon,1},\ldots,g^\I_{\varepsilon,n})\in\mathbb{R}^n$ by
 \begin{equation}\label{Ideal Direction}
 g^\I_{\varepsilon,i}:= \text{argmin}\left\{ \lvert g \rvert \,\,\, : \,\,\, g\in G^m_{\varepsilon,i}({\mathbf x})  \right\}, \qquad i=1,\dots,n.
 \end{equation} 
 In fact, the $i$-th component of $\mathbf g^\I_\varepsilon$ denotes the member of $G^m_{\varepsilon,i}({\mathbf x})$ with minimum distance from zero. It is easy to see that $g^\I_{\varepsilon,i}$ can be expressed by
$$ g^\I_{\varepsilon,i}= \frac{{\rm sign}(m_i)+{{\rm sign}}(M_i)}{2} \min\left\{\lvert m_i \rvert , \lvert M_i\rvert \right\}, \qquad i=1,\ldots,n, $$
where
$$ m_i := \min\left\{\frac{\partial f}{\partial x_i}(\mathbf s_j) \,\, : \,\, j=0,\ldots,m \right\} ,$$
and
$$ M_i:=\max\left\{\frac{\partial f}{\partial x_i}(\mathbf s_j) \,\, : \,\, j=0,\ldots,m \right\}. $$
Therefore, in order to obtain the $i$-th component of $\mathbf g^\I_\varepsilon$, we only need to compute the partial derivative of  function $f$ with respect to $x_i$ at sampled points $\mathbf s_0,\ldots,\mathbf s_m$ and find the ones with minimum and maximum value. Thus, there is no need to consider any complex problem, like QP \eqref{CQP}, in order to obtain $\mathbf g^\I_\varepsilon$ and hence these directions are easy to compute. The preceding discussion is summarized in the following definition.
\begin{definition}\label{Ideal_D}
	Let $\mathbf g^\I_\varepsilon$ be as defined in \eqref{Ideal Direction}. If $\mathbf g^\I_\varepsilon\neq \mathbf 0$, the direction $\mathbf d^\I_\varepsilon=(d^\I_{\varepsilon,1},\ldots,d^\I_{\varepsilon,n}):=- \lVert \mathbf g^\I_\varepsilon\rVert^{-1}\mathbf g^\I_\varepsilon$ is called (normalized) Ideal direction of the set $G^m_\varepsilon(\mathbf x)$.   
\end{definition}
In what follows we examine the most important properties of Ideal directions.
\begin{lemma}\label{lemma1}
	Let ${\mathbf x}$ be a differentiable point, $m\in\mathbb{N}$ and $\varepsilon>0$.
	\begin{itemize}
	\item[1.]For the set $G^m_\varepsilon({\mathbf x})$ we always have
		$$\lVert {\mathbf g}^\I_\varepsilon\rVert\leq\lVert {\mathbf g}^\s_\varepsilon\rVert. $$
	\item[2.] If $\mathbf 0\in G^m_\varepsilon({\mathbf x})$, then
		$${\mathbf g}^\I_\varepsilon= \mathbf 0 .$$   
		\end{itemize}
\end{lemma}
{\it Proof}
	\begin{itemize}
		\item[1.] Since ${\mathbf g}^\s_\varepsilon$ is the solution of the QP \eqref{CQP}, we have ${\mathbf g}^\s_\varepsilon\in G^m_\varepsilon({\mathbf x})$ and hence $ g^\s_{\varepsilon,i}\in G^m_{\varepsilon,i}({\mathbf x})$ for every $i=1,\ldots,n$. Now, since $ g^\I_{\varepsilon,i}$ is the solution of problem \eqref{Ideal Direction} and $ g^\s_{\varepsilon,i}\in G^m_{\varepsilon,i}({\mathbf x})$, we conclude that
		$$\lvert  g^\I_{\varepsilon,i}\rvert\leq \lvert  g^\s_{\varepsilon,i}\rvert, \qquad \text{for \,all} \,\, i=1,\ldots,n,$$
		yielding $\lVert {\mathbf g}^\I_\varepsilon\rVert\leq\lVert {\mathbf g}^\s_\varepsilon\rVert.$
		\item[2.] Since $\mathbf 0\in G^m_\varepsilon({\mathbf x})$, we conclude that ${\mathbf g}^\s_\varepsilon=\mathbf 0$ and the result follows immediately from the first part of this lemma. \qed
	\end{itemize}

In the following lemma, we show that if ${\mathbf g}^\I_\varepsilon\neq \mathbf 0$, the Ideal direction ${\mathbf d}^\I_\varepsilon=-\lVert \mathbf g^\I_\varepsilon\rVert^{-1}\mathbf g^\I_\varepsilon$ is an ADD.
\begin{lemma}\label{Armijo_C}
	Suppose that ${\mathbf g}^\I_\varepsilon\neq \mathbf 0$ for the set $G^m_\varepsilon({\mathbf x})$ and $c\in (0,1)$. Then there exists $z_0>0$ such that
	$$ f({\mathbf x}+t\mathbf d^\I_\varepsilon)-f({\mathbf x}) < -tc\lVert {\mathbf g}^\I_\varepsilon\rVert , \qquad {\rm{for \,\,all}} \,\, t\in(0,z_0].$$
\end{lemma}
{\it Proof}
	Let $ A_i:= G^m_{\varepsilon,i}({\mathbf x})$. By noting that $\frac{\partial f}{\partial x_i}(\mathbf x)\in A_i$ for every $i\in\{1,\ldots,n\}$, we have
	$$ \nabla f({\mathbf x})^t {\mathbf d}^\I_\varepsilon=\sum_{i=1}^n \frac{\partial f}{\partial x_i}({\mathbf x}) d^\I_{\varepsilon,i} \leq \sum_{i=1}^n \sup_{g\in A_i} \{g\, d^\I_{\varepsilon,i}\} =-\lVert {\mathbf g}^\I_\varepsilon \rVert^{-1} \sum_{i=1}^n \sup_{g\in A_i} \{g\,  g^\I_{\varepsilon,i}\}. $$
	Since $g^\I_{\varepsilon,i}$ is the solution of the problem \eqref{Ideal Direction}, a necessary condition for minimum implies that for every $i=1,\ldots,n$ 
	$$ g\,  g^\I_{\varepsilon,i}  \geq  ({g^\I_{\varepsilon,i}})^2 , \qquad \text{for \,\,all} \,\, g\in A_i,$$
	therefore 
	$$ \nabla f({\mathbf x})^t {\mathbf d}^\I_\varepsilon \leq -\lVert {\mathbf g}^\I_\varepsilon \rVert^{-1} \sum_{i=1}^n ({g^\I_{\varepsilon,i}})^2 =-\lVert {\mathbf g}^\I_\varepsilon \rVert,  $$
	and the result follows from the fact that $\nabla f({\mathbf x})^t {\mathbf d}^\I_\varepsilon=f'({\mathbf x},{\mathbf d}^\I_\varepsilon)$.\qed

Therefore, when ${\mathbf g}^\I_\varepsilon\neq \mathbf 0$ an ADD is in hand without solving any complex subproblem. We will use this key property to alleviate the first drawback we mentioned above. In the next lemma, let $\mathbf e_i$ be a vector with a $1$ in the $i$-th coordinate and zeros elsewhere.

\begin{lemma}
	For a differentiable point ${\mathbf x}\in\mathbb{R}^n, \varepsilon>0$ and $m\in\mathbb{N}$, let ${\mathbf d}^\I_\varepsilon$ be the Ideal direction of the set $G^m_\varepsilon({\mathbf x})$. Then ${\mathbf d}^\I_\varepsilon\neq \mathbf 0$  if and only if, for some $i\in\{1,\ldots,n\}$ there exists a hyperplane $H$ with the normal vector $\mathbf e_i$ which  strongly separates $G^m_\varepsilon({\mathbf x})$ and origin.  
\end{lemma}
{\it Proof}
	First, assume that for some $i\in\{1,\dots,n\}$ there is a hyperplane $H$ with the normal vector $\mathbf e_i$ strongly separating $G^m_\varepsilon({\mathbf x})$ and origin. This means that $0\notin G^m_{\varepsilon,i}({\mathbf x})$ and hence
	$$ g^\I_{\varepsilon,i}=\text{argmin}\left\{ \lvert g \rvert \,\,\, : \,\,\, g\in G^m_{\varepsilon,i}({\mathbf x})  \right\}\neq 0,$$
	which yields ${\mathbf g}^\I_\varepsilon\neq \mathbf 0$ and  consequently ${\mathbf d}^\I_\varepsilon\neq \mathbf 0$. Now, suppose that ${\mathbf d}^\I_\varepsilon\neq \mathbf 0$ for the set $G^m_\varepsilon({\mathbf x})$. Thus, ${\mathbf g}^\I_\varepsilon\neq \mathbf 0$ and there exits $i\in\{1,\ldots,n\}$ such that $ g^\I_{\varepsilon,i}\neq0$ which means that $0\notin G^m_{\varepsilon,i}({\mathbf x})$. Therefore, by convexity of $G^m_\varepsilon({\mathbf x})$, we conclude that the set $G^m_\varepsilon({\mathbf x})$ does not touch the hyperplane
	$$H=\{(h_1,\ldots,h_{i-1}, g^\I_{\varepsilon,i}/2,h_{i+1},\dots,h_n) \,\, : \,\, h_j\in\mathbb{R}, \,\,\, \text{for \,all}\,\, j\neq i\}.$$
	Clearly, $\mathbf e_i$ is the normal vector of $H$ and this hyperplane strongly separates $G^m_\varepsilon({\mathbf x})$ and origin.\qed

The following lemma states that, if function $f$ is smooth at a nonstationary  point ${\mathbf x}\in\mathbb{R}^n$, there exists a sampling radius $\bar{\varepsilon}>0$ such that the Ideal direction of the set $G^m_{\bar{\varepsilon}}({\mathbf x})$ is nonzero.
\begin{lemma}\label{Lemma2}
	Suppose that the function $f:\mathbb{R}^n\to\mathbb{R}$ is smooth at a point ${\mathbf x}\in\mathbb{R}^n$ and $\nabla f({\mathbf x})\neq \mathbf 0$. Then, there exists a sampling radius $\bar{\varepsilon}>0$ such that ${\mathbf d}^\I_{\bar{\varepsilon}}\neq \mathbf 0$ for the set $G^m_{\bar{\varepsilon}}({\mathbf x})$.
\end{lemma}
{\it Proof}
	Since the map $\nabla f:D\subset\mathbb{R}^n\to\mathbb{R}^n$ is continuous at the point ${\mathbf x}\in\mathbb{R}^n$, for every $\delta>0$ there exists $\varepsilon>0$ such that
	$$\lVert \nabla f({\mathbf x})-\nabla f(\mathbf y)\rVert < \delta, \qquad \text{for \,all} \,\,\mathbf y\in B({\mathbf x},\varepsilon). $$
	Now since $\nabla f({\mathbf x})\neq \mathbf 0$, without loss of generality, there exists $i\in \{1,\ldots,n\}$ such that $\frac{\partial f({\mathbf x})}{\partial x_i}>0$. By taking $\bar{\delta}:=\frac{\partial f({\mathbf x})}{\partial x_i}/2$, there exists $\bar{\varepsilon}>0$ such that
	\begin{equation}\label{Ideal effieciency}
	\left\lvert \frac{\partial f({\mathbf x})}{\partial x_i}-\frac{\partial f(\mathbf y)}{\partial y_i}\right\lvert \leq \lVert \nabla f({\mathbf x})-\nabla f(\mathbf y)\rVert< \frac{1}{2} \frac{\partial f({\mathbf x})}{\partial x_i}, \qquad \text{for \,all}\, \, \mathbf y\in B({\mathbf x},\bar{\varepsilon}).
	\end{equation}
	This means that $\frac{\partial f(\mathbf y)}{\partial y_i}>0$ for every $\mathbf y\in B({\mathbf x},\bar{\varepsilon})$ and consequently $0\notin G^m_{\bar{\varepsilon},i}({\mathbf x})$ that gives  $ g^\I_{\bar{\varepsilon},i}\neq 0$. Therefore ${\mathbf d}^\I_{\bar{\varepsilon}}=-\lVert \mathbf g^\I_{\bar{\varepsilon}}\rVert^{-1}{\mathbf g}^\I_{\bar{\varepsilon}}\neq \mathbf 0$  for the set $G^m_{\bar{\varepsilon}}({\mathbf x})$. \qed

In the light of preceding lemma, one might think that, whenever we obtain a zero Ideal direction, then we can reduce the sampling radius to  obtain a nonzero one. However, reducing sampling radius may lead to losing information of the nearby nonsmooth curve and hence the resulting nonzero Ideal direction would not be an effective search direction. In such cases, the approximate $\varepsilon$-steepest descent direction obtained by QP \eqref{CQP} may provide a more suitable search direction. 

By the proof of the Lemma \ref{Lemma2} one can see that, when we are far away from a nonsmooth curve, in the sense that $B({\mathbf x},\varepsilon)\subset D$, since the information we collect by sampling gradients are very close together, we expect the inequalities \eqref{Ideal effieciency} hold (without reducing the sampling radius) and for the set $G^m_\varepsilon({\mathbf x})$ we obtain a nonzero Ideal direction. This important observation resolves the second drawback stated above.

It is noted that, the Ideal direction is simply obtained  without solving any time-consuming subproblem. In fact, computing
gradient at sampled points is equivalent to having an Ideal direction and whenever ${\mathbf g}^\I_\varepsilon\neq \mathbf 0$ it decreases the value of function along ${\mathbf d}^\I_\varepsilon$ with a step size $t$ satisfying the Armijo condition.
In the case that ${\mathbf g}^\I_\varepsilon=\mathbf 0$ for the set $ G^m_\varepsilon ({\mathbf x})$, according to the second part of Lemma \ref{lemma1},   we consider ${\mathbf x}$ as a candidate to be a Clarke $\varepsilon$-stationary point and this possibility is checked by solving the QP \eqref{CQP}.

\section{A Gradient Sampling method using Ideal direction }\label{Section5}	

In this section, based on Ideal directions, we modify the GS algorithm presented in \cite{Kiwiel2007}. Next, following the work of Kiwiel in \cite{Kiwiel2007}, we analyze the global  convergence of the proposed method.

To modify the GS method, our method consider the Ideal direction as the first choice for the search direction. As mentioned, whenever ${\mathbf g}^\I_\varepsilon\neq \mathbf 0$ it presents an ADD without solving any complex subproblem. Thus, it is reasonable to compute the Ideal direction before solving the QP \eqref{CQP}. Furthermore, during each iteration, the optimality condition \eqref{Opt_Tol} must be checked. By noting that $\lVert {\mathbf g}^\I_{\varepsilon_k} \rVert\leq \lVert {\mathbf g}^\s_{\varepsilon_k} \rVert$, if for a point ${\mathbf x}_k\in\mathbb{R}^n$ we have  $\lVert {\mathbf g}^\I_{\varepsilon_k} \rVert >\nu_k$, we conclude that $\lVert {\mathbf g}^\s_{\varepsilon_k}\lVert>\nu$ and hence ${\mathbf x}_k$ is not a $(\nu_k , G^m_{\varepsilon_k}({\mathbf x}_k))$-stationary point. In the case that $\lVert {\mathbf g}^\I_{\varepsilon_k}\lVert\leq\nu_k$, we need to check whether ${\mathbf x}_k$ is $\left(\nu_k , G^m_{\varepsilon_k}({\mathbf x}_k)\right)$-stationary point or not and this is done by solving the QP \eqref{CQP}. Therefore, in some cases we can check the optimality condition \eqref{Opt_Tol} without computing ${\mathbf g}^\s_{\varepsilon_k}$. In this way, through Ideal directions, we can reduce the number of quadratic subproblems significantly. 

Based on the preceding discussion, a Gradient Sampling algorithm based on Ideal directions (GSI) is presented in Algorithm \ref{Algorithm 2} .

\IncMargin{1em}
\begin{algorithm}

	\SetKwInOut{Input}{inputs}
	\SetKwInOut{Output}{output}
	\SetKwInOut{Required}{requied}
	
	\Indm  
	\Input {${\mathbf x}_0\in\mathbb{R}^n$ as initial guess, 
      			initial sampling radius $\varepsilon_0>0$, initial stationarity tolerance $\nu_0> 0$, sample size $m\in\mathbb{N}$, 
				parameters $\nu_{opt}, \varepsilon_{opt}\geq0$ as tolerances in stopping condition, 
				parameters $\theta,\mu\in(0,1)$ as reduction factors for optimality tolerance and sampling radius, 
				and the backtracking Armijo line search parameters $\gamma, c\in(0,1)$.

	}
	\Output{An approximation of a Clarke stationary point of $f$.}
	\Required{$f$ is differentiable at the starting point ${\mathbf x}_0\in\mathbb{R}^n$.}
	\Indp
	\BlankLine
	$k:=0$\;
	\While{$\nu_k\geq\nu_{opt}$ {\rm or} $\varepsilon_k\geq \varepsilon_{opt}$}{ 
		
						Sample $\mathbf u_{k1},\ldots,\mathbf u_{km}$  independently and uniformly from $B_1$ and set
						$ \mathbf s_{kj}:={\mathbf x}_k+\varepsilon_k \mathbf u_{kj}$ for $j=1,\ldots,m$\;

						 \If{ $\mathbf s_{kj}$ is not a differentiable point for some $j\in\{1,\ldots,m\}$}
						 {
						 	\textbf{Stop}\; 
					 	}
						 Compute $\mathbf g_{\varepsilon_k}^\I$ for the set $G^m_{\varepsilon_k}({\mathbf x}_k)$ and set $\mathbf g_k:=\mathbf g_{\varepsilon_k}^\I$\;
						 
						 \If{$\lVert \mathbf g_k \rVert \le \nu_k$}
						 {
						 		Solve QP \eqref{CQP} to compute $\mathbf g^\s_{\varepsilon_k}$  and set $\mathbf g_k:=\mathbf g^\s_{\varepsilon_k} $\;

						 }
						 \If{$\lVert \mathbf g_k\rVert\leq\nu_{k}$}
						 {$\nu_{k+1}:=\theta \nu_{k}$, $\varepsilon_{k+1}:=\mu \varepsilon_{k}$\;
						 	$\mathbf x_{k+1}:=\mathbf x_k$\;
						 	$k:=k+1$\;
						 	\textbf{Continue}\;
						 	
						 	}
						  	$\mathbf d_k:=-\lVert \mathbf g_k\rVert^{-1} \mathbf g_k$\;
						 	$t_k:=\texttt{BALS}({\mathbf x}_k, \mathbf g_k, \mathbf d_k,\gamma,c)$\;
						 	${\mathbf x}_{k+1}:={\mathbf x}_k+t_k \mathbf d_k $\;

						 \If{$f$ is not differentiable at ${\mathbf x}_{k+1}$}
						 {
						 		Find a point $\hat{\mathbf x}$ at which $f$ is differentiable and satisfies
						 		
						 	   $\qquad  f(\hat{{\mathbf x}})-f({\mathbf x}_k)<-ct_k\lVert \mathbf g_k\rVert$\; 
						 		$\qquad  \lVert {\mathbf x}_k+t_k  \mathbf d_k-\hat{{\mathbf x}}\rVert \leq \min\{t_k,\varepsilon_k\}$\;
						 		Set  $\mathbf x_{k+1}:=\hat{\mathbf x}$\;
					 	 }
				 	    $ \nu_{k+1}:=\nu_k, \varepsilon_{k+1}:=\varepsilon_k$\;		  
						$k:=k+1$\;

	}
   Introduce ${\mathbf x}_k$ as an approximation of a Clarke stationary point.
	\caption{Gradient Sampling method using Ideal direction (GSI) }
	\label{Algorithm 2}
\end{algorithm}
\DecMargin{1em}

\subsection{Convergence analysis}

To study the convergence of the GSI algorithm, at first we recall three notations from \cite{Kiwiel2007,Burke2005}.

 For $\varepsilon>0$ and $\bar{{\mathbf x}}\in\mathbb{R}^n$, let $\rho_{\varepsilon}(\bar{{\mathbf x}})$ be the distance between $G_{\varepsilon}(\mathbf{\bar{x}})$ and origin. In mathematical terms
$$\rho_{\varepsilon}(\bar{{\mathbf x}}):= \text {dist}\Big(\mathbf 0,G_\varepsilon(\bar{{\mathbf x}})\Big), $$
and for ${\mathbf x}\in\mathbb{R}^n$ and $m\in\mathbb{N}$, let
$$D^m_\varepsilon({\mathbf x}):=\prod_1^m \left(B({\mathbf x},\varepsilon)\cap D\right)\subset \prod_1^m \mathbb{R}^n. $$
Moreover, for $\mathbf x, \bar{\mathbf x}\in\mathbb{R}^n, m\in\mathbb{N}$ and $\varepsilon, \delta>0$ define

$$V^m_\varepsilon(\bar{{\mathbf x}},{\mathbf x},\delta):=\{(\mathbf y_1,\ldots,\mathbf y_m)\in D^m_\varepsilon({\mathbf x}) \,: \, \text{dist}(\mathbf 0,\texttt{co}\{\nabla f(\mathbf y_k)\}_{k=1}^m)\leq \rho_\varepsilon(\bar{{\mathbf x}})+\delta \}. $$

Next, we recall the next three lemmas from \cite{Kiwiel2007}. Before it, suppose that the objective function in problem \eqref{Main-Problem} fulfills Assumption \ref{Assumption} throughout this section \cite{GS-full}.
   
   \begin{assumption}\label{Assumption}
   	The function $f:\mathbb{R}^n\to\mathbb{R}$ is locally Lipschitz and continuously differentiable on an open set $D$ with full measure in $\mathbb{R}^n$.
   \end{assumption}
In the next three lemmas, the first one provides a variational inequality for some elements of $G_\varepsilon(\bar{{\mathbf x}})$ and the second one  examines the local behavior of the set $V^m_\varepsilon(\bar{{\mathbf x}},\cdot,\delta)$ in the vicinity of $\bar{{\mathbf x}}$. Finally, Lemma \ref{K3}  derives a lower bound for the step size $t_k$ when the sampled points $(\mathbf s_{k1},\ldots,\mathbf s_{km})$ lie in a special open subset of $V^m_{\varepsilon_k}(\bar{{\mathbf x}},{\mathbf x}_k,\delta)$. 

\begin{lemma}[\cite{Kiwiel2007}]\label{K1}
	For a point $\bar{{\mathbf x}}\in \mathbb{R}^n$,	suppose that $\mathbf 0\notin G_{\varepsilon}(\bar{{\mathbf x}})$ and $c\in (0,1)$, then there is $\delta>0$ such that for every $\mathbf u\in G_\varepsilon(\bar{{\mathbf x}})$ satisfying
	$\lVert \mathbf u\rVert\leq \rho_\varepsilon(\bar{{\mathbf x}})+\delta, $
	we have
	$$\langle \mathbf u,\mathbf v\rangle > c \lVert \mathbf u \rVert^2, \qquad {\rm{for \,\, all}} \,\, \mathbf v\in G_\varepsilon(\bar{{\mathbf x}}).$$ 
\end{lemma}

\begin{lemma}[\cite{Kiwiel2007}]\label{K2}
	For any $\delta>0$, there exist $\tau>0, \overline{m}\geq n+1$ and a nonempty open set $\overline{V}$, such that 
	$$\overline{V}\subset V^{\overline{m}}_\varepsilon(\bar{{\mathbf x}},{\mathbf x},\delta), \qquad {\rm{for\, all}} \,\, {\mathbf x}\in B(\bar{{\mathbf x}},\tau).$$

\end{lemma}

\begin{lemma}[\cite{Kiwiel2007}]\label{K3}
	Assume that $\mathbf 0\notin G_\varepsilon(\bar{{\mathbf x}}) $, choose $\delta>0$ as in Lemma \ref{K1} and $\tau, \overline{m}, \overline{V}$ as in Lemma \ref{K2}. Suppose that at iteration $k$ of Algorithm \ref{Algorithm 2}, $(\mathbf s_{k1},\ldots,\mathbf s_{k\overline{m}})\in \overline{V}$ and the backtracking Armijo line search in Line 18 is reached with ${\mathbf x}_k\in B(\bar{{\mathbf x}},\min\{\tau,\varepsilon_k/3\})$ and $\mathbf d_k=\mathbf d^\s_{\varepsilon_k}$. Then $t_k\geq \min\{1,\gamma \varepsilon_k/3\}$.  
\end{lemma}

In order to extend Lemma \ref{K3} to Ideal directions, at first we need to introduce the following notations. For $\mathbf{\bar{x}}\in\mathbb{R}^n$ and $i\in\{1,\ldots,n\}$, let
$$\rho_{\varepsilon,i}(\bar{{\mathbf x}}):= \text {dist}\Big(0,G_{\varepsilon,i}(\bar{{\mathbf x}})\Big), $$
where
$$G_{\varepsilon,i}(\bar{{\mathbf x}}):= \texttt{cl\;co}\left\{\frac{\partial f}{\partial x_i}(B(\bar{{\mathbf x}},\varepsilon)\cap D)\right\}.$$ 
Furthermore, for $\mathbf{x,\bar{x}}\in\mathbb{R}^n, m\in\mathbb{N}$ and $\varepsilon, \delta,\omega>0$ we define

\begin{align*}
&W^m_\varepsilon(\bar{{\mathbf x}},{\mathbf x},\delta,\omega):=\\
&\ \left\{
(\mathbf y_1,\ldots,\mathbf y_m)\in D^m_\varepsilon({\mathbf x}) \, : \, 
\begin{aligned}
& \text{dist}(0,\texttt{co}\{\tfrac{\partial f}{\partial x_i}(\mathbf y_k)\}_{k=1}^m)\leq \rho_{\varepsilon,i}(\bar{{\mathbf x}})+\delta, \,\,  i\in\mathcal{A}(\bar{{\mathbf x}}) \text{ and} \\ & 
\text{dist}(0,\texttt{co}\{\tfrac{\partial f}{\partial x_i}(\mathbf y_k)\}_{k=1}^m)=0,  \,\,  i\in\mathcal{A}'_\texttt{int}(\bar{{\mathbf x}})\text{ and} \\ &
\text{dist}(0,\texttt{co}\{\tfrac{\partial f}{\partial x_i}(\mathbf y_k)\}_{k=1}^m)\leq\omega,  \,\,  i\in\mathcal{A}'_\texttt{bd}(\bar{{\mathbf x}})
\end{aligned}
\right\},
\end{align*}
in which 
\begin{align*}
&\mathcal{A}(\bar{{\mathbf x}}):=\{i\in\{1,\ldots,n\} \, : \, 0\notin G_{\varepsilon,i}(\bar{{\mathbf x}})\},\\&
\mathcal{A}'_\texttt{int}(\bar{{\mathbf x}}):=\{i\in\{1,\ldots,n\} \, : \, 0\in \texttt{int} \, G_{\varepsilon,i}(\bar{{\mathbf x}})\},\\&
\mathcal{A}'_\texttt{bd}(\bar{{\mathbf x}}):=\{i\in\{1,\ldots,n\} \, : \, 0\in \texttt{bd} \, G_{\varepsilon,i}(\bar{{\mathbf x}})\}.
\end{align*}
Here $\texttt{bd}$ and $\texttt{int}$ denote the boundary and interior of a set respectively.

To study the global convergence behavior of our method, we have to understand the local behavior of the set $W^m_\varepsilon(\bar{{\mathbf x}},\cdot,\delta,\omega)$ in the vicinity of $\bar{{\mathbf x}}$. Then, we derive a lower bound for the step size $t_k$ when the search direction is defined by an Ideal direction and the sampled points $(\mathbf s_{k1},\ldots,\mathbf s_{km})$ are located in a special open subset of $W^m_{\varepsilon_k}(\bar{{\mathbf x}},{\mathbf x}_k,\delta,\omega)$. To this end, we start with the following lemma which is an immediate consequence of Lemma \ref{K1}.

\begin{lemma}\label{zero p}
	For a point $\bar{{\mathbf x}}\in \mathbb{R}^n, c\in (0,1)$ and any $i\in\mathcal{A}(\bar{{\mathbf x}})$, there is $\delta_i>0$ such that for every $u_i\in G_{\varepsilon,i}(\bar{{\mathbf x}})$ satisfying
	$\lvert u_i\rvert\leq \rho_{\varepsilon,i}(\bar{{\mathbf x}})+\delta_i, $
	we have
	$$ u_i v_i > c  u_i^2, \qquad {\rm{for \,\, all}} \,\, v_i\in G_{\varepsilon,i}(\bar{{\mathbf x}}).$$ 
\end{lemma}
{\it Proof}
	This follows immediately from Lemma \ref{K1}.\qed

The following lemma states that for the points sufficiently close to $\bar{{\mathbf x}}$, the set valued map $W^m_\varepsilon(\bar{{\mathbf x}},\cdot,\delta,\omega)$ contains a nonempty open set.
\begin{lemma}\label{first p}
	For any $\varepsilon, \delta, \omega>0$, there exist $\tau>0, \overline{m}\in\mathbb{N}$ and a nonempty open set $\overline{W}$, such that 
	 $$\overline{W}\subset W^{\overline{m}}_\varepsilon(\bar{{\mathbf x}},{\mathbf x},\delta,\omega),\qquad {\rm{for \,\,all}} \,\, {\mathbf x}\in B(\bar{{\mathbf x}},\tau).$$
\end{lemma}
{\it Proof}\, First, assume that $i\in\mathcal{A}(\bar{{\mathbf x}})$. Since $\delta>0$, there exists $u_i\in \texttt{co}\left\{\frac{\partial f}{\partial x_i}(B(\bar{{\mathbf x}},\varepsilon)\cap D)\right\}$ such that
	$$\lvert u_i\rvert \leq \rho_{\varepsilon,i}(\bar{{\mathbf x}})+\delta , \qquad \text{for\, all}\,\, i\in\mathcal{A}(\bar{{\mathbf x}}), $$
	then, by Carath\'{e}odory's theorem, there exist $\mathbf s_i^1, \mathbf s_i^2\in B(\bar{{\mathbf x}},\varepsilon)\cap D$ and nonnegative scalars $\bar{\lambda}_i^1, \bar{\lambda}_i^2$ such that $\bar{\lambda}_i^1+\bar{\lambda}_i^2=1$ and
	$$u_i=\bar{\lambda}_i^1\frac{\partial f}{\partial x_i}(\mathbf s_i^1)+ \bar{\lambda}_i^2\frac{\partial f}{\partial x_i}(\mathbf s_i^2). $$
	Now according to Assumption \ref{Assumption}, for each $i\in\{1,\ldots,n\}$ the map $\frac{\partial f}{\partial x_i}(\cdot)$ is continuous in $D$. Thus, for every $i\in\mathcal{A}(\bar{{\mathbf x}})$ there is $\bar{\varepsilon}_1\in(0,\varepsilon)$ such that for all $\mathbf y_i^1\in B(\mathbf s_i^1,\bar{\varepsilon}_1)\cap D$ and $\mathbf y_i^2\in B(\mathbf s_i^2,\bar{\varepsilon}_1)\cap D$, we have
	\begin{equation}\label{c1}
	\Big\lvert \bar{\lambda}_i^1\frac{\partial f}{\partial x_i}(\mathbf y_i^1)+ \bar{\lambda}_i^2\frac{\partial f}{\partial x_i}(\mathbf y_i^2)\Big\rvert\leq \rho_{\varepsilon,i}(\bar{{\mathbf x}})+\delta, \qquad \text{for\, all}\,\, i\in\mathcal{A}(\bar{{\mathbf x}}). 
	\end{equation}
	Second, suppose that $i\in \mathcal{A}'_\texttt{int}(\bar{{\mathbf x}})$. Then we have $0\in\texttt{co}\left\{ \frac{\partial f}{\partial x_i}(B(\bar{{\mathbf x}},\varepsilon)\cap D)\right\}$ which implies the existence of $0<u_i^+\in \texttt{co}\left\{ \frac{\partial f}{\partial x_i}(B(\bar{{\mathbf x}},\varepsilon)\cap D)\right\}$ and $0>u_i^-\in \texttt{co}\left\{ \frac{\partial f}{\partial x_i}(B(\bar{{\mathbf x}},\varepsilon)\cap D)\right\}$. Using Carath\'{e}odory's theorem, for each $i\in\mathcal{A}'_\texttt{int}(\bar{{\mathbf x}})$, there exist $\mathbf s_i^1, \mathbf s_i^2, \mathbf s_i^3, \mathbf s_i^4\in B(\bar{{\mathbf x}},\varepsilon)\cap D$ and nonnegative scalars $\bar{\lambda}_i^1, \bar{\lambda}_i^2, \bar{\lambda}_i^3, \bar{\lambda}_i^4$ such that
	$$0<u_i^+= \bar{\lambda}_i^1\frac{\partial f}{\partial x_i}(\mathbf s_i^1)+ \bar{\lambda}_i^2\frac{\partial f}{\partial x_i}(\mathbf s_i^2), \qquad \bar{\lambda}_i^1+\bar{\lambda}_i^2=1,  $$
	and
	$$0>u_i^-= \bar{\lambda}_i^3\frac{\partial f}{\partial x_i}(\mathbf s_i^3)+ \bar{\lambda}_i^4\frac{\partial f}{\partial x_i}(\mathbf s_i^4), \qquad \bar{\lambda}_i^3+\bar{\lambda}_i^4=1. $$
	Now since the map $\frac{\partial f}{\partial x_i}(\cdot)$ is continuous in $D$, there is $\bar{\varepsilon}_2\in(0,\varepsilon)$ such that for every $i\in\mathcal{A}'_\texttt{int}(\bar{{\mathbf x}})$ one can write
	\begin{equation}\label{c2}
	0<\bar{\lambda}_i^1\frac{\partial f}{\partial x_i}(\mathbf y_i^1)+\bar{\lambda}_i^2\frac{\partial f}{\partial x_i}(\mathbf y_i^2), \qquad \text{for\, all }\, \mathbf y_i^1\in B(\mathbf s_i^1,\bar{\varepsilon}_2), \, \mathbf y_i^2\in B(\mathbf s_i^2,\bar{\varepsilon}_2),
	\end{equation}
	and
	\begin{equation}\label{c3}
	0>\bar{\lambda}_i^3\frac{\partial f}{\partial x_i}(\mathbf y_i^3)+\bar{\lambda}_i^4\frac{\partial f}{\partial x_i}(\mathbf y_i^4), \qquad \text{for\, all }\, \mathbf y_i^3\in B(\mathbf s_i^3,\bar{\varepsilon}_2), \, \mathbf y_i^4\in B(\mathbf s_i^4,\bar{\varepsilon}_2).
	\end{equation}
	Finally, let $i\in\mathcal{A}'_\texttt{bd}(\bar{{\mathbf x}})$. Since $\omega>0$, there is $u_i\in \texttt{co}\left\{ \frac{\partial f}{\partial x_i}(B(\bar{{\mathbf x}},\varepsilon)\cap D)\right\}$ such that
	$$\lvert u_i\rvert \leq \rho_{\varepsilon,i}(\bar{{\mathbf x}})+\omega=\omega , \qquad \text{for\, all}\,\, i\in\mathcal{A}'_\texttt{bd}(\bar{{\mathbf x}}). $$
	By Carath\'{e}odory's theorem, for each $i\in\mathcal{A}'_\texttt{bd}(\bar{{\mathbf x}})$, there exist $\mathbf s_i^1, \mathbf s_i^2\in B(\bar{{\mathbf x}},\varepsilon)\cap D$ and nonnegative scalars $\bar{\lambda}_i^1, \bar{\lambda}_i^2$ such that $\bar{\lambda}_i^1+\bar{\lambda}_i^2=1$ and
	$$u_i=\bar{\lambda}_i^1\frac{\partial f}{\partial x_i}(\mathbf s_i^1)+ \bar{\lambda}_i^2\frac{\partial f}{\partial x_i}(\mathbf s_i^2).$$
	Using continuity of the map $\frac{\partial f}{\partial x_i}(\cdot)$, there is $\bar{\varepsilon}_3\in(0,\varepsilon)$ such that for all $\mathbf y_i^1\in B(\mathbf s_i^1,\bar{\varepsilon}_3)$ and $\mathbf y_i^2\in B(\mathbf s_i^2,\bar{\varepsilon}_3)$ we have
	\begin{equation}\label{c4}
	\Big\lvert \bar{\lambda}_i^1\frac{\partial f}{\partial x_i}(\mathbf y_i^1)+ \bar{\lambda}_i^2\frac{\partial f}{\partial x_i}(\mathbf y_i^2)\Big\rvert\leq \omega.
	\end{equation}
	Now there exists $0<\bar{\varepsilon}\leq \min \{\bar{\varepsilon}_1,\bar{\varepsilon}_2,\bar{\varepsilon}_3\}$ such that $\overline{W}:=\prod_{i,j}\texttt{int}B(\mathbf s_i^j,\bar{\varepsilon}) $ lies in $D^{\overline{m}}_{\varepsilon-\bar{\varepsilon}}(\bar{{\mathbf x}})$ in which 
	$$2n \leq\overline{m}:=2 (\lvert \mathcal{A}(\bar{{\mathbf x}})\rvert+\lvert\mathcal{A}'_\texttt{bd}(\bar{{\mathbf x}})\rvert) +4 \lvert \mathcal{A}'_\texttt{int}(\bar{{\mathbf x}})\rvert \leq 4n. $$
	Furthermore, in view of \eqref{c1}, \eqref{c2}, \eqref{c3} and \eqref{c4} and the fact that $B(\bar{{\mathbf x}},\varepsilon-\tau)\subset B({\mathbf x},\varepsilon)$ for all ${\mathbf x}\in B(\bar{{\mathbf x}},\tau)$ with $\tau:=\bar{\varepsilon}$, we can conclude that
	$$\overline{W}\subset W^{\overline{m}}_\varepsilon(\bar{{\mathbf x}},{\mathbf x},\delta,\omega), \qquad \text{for\, all}\,\, {\mathbf x}\in B(\bar{{\mathbf x}},\tau).$$ \qed

The following corollary is an immediate consequence of the preceding lemma. 

\begin{corollary}\label{corollary zero}
	Let $\varepsilon, \delta, \omega >0$, $\bar{{\mathbf x}},{\mathbf x}\in\mathbb{R}^n$ and $\overline{W}\subset W^{\overline{m}}_{\varepsilon}(\bar{{\mathbf x}},{\mathbf x},\delta,\omega)$ be as in Lemma \ref{first p}. For the set $G^{\overline{m}}_\varepsilon(\mathbf x)$ suppose that $(\mathbf s_1,\ldots,\mathbf s_{\overline{m}})\in\overline{W}$ and $ \mathbf g^\I_\varepsilon=(g^\I_{\varepsilon,1},\ldots,g^\I_{\varepsilon,n})$ is defined as \eqref{Ideal Direction}. Then

	\begin{align*}
	& \lvert g^\I_{\varepsilon,i}\rvert\leq\rho_{\varepsilon,i}(\bar{{\mathbf x}})+\delta, \qquad {\rm{for\, all}} \,\, i\in\mathcal{A}(\bar{{\mathbf x}}), \\
	& g^\I_{\varepsilon,i}=0, \qquad {\rm{for\, all}} \,\, i\in\mathcal{A}'_{\rm{{int}}}(\bar{{\mathbf x}}), \\
	& \lvert g^\I_{\varepsilon,i}\rvert\leq \omega, \qquad {\rm{for\, all}} \,\, i\in\mathcal{A}'_{\rm{{bd}}}(\bar{{\mathbf x}}).
	\end{align*}
\end{corollary}

The next lemma indicates that, if the set of sampled points lies in the open set $\overline{W}$ and the search direction is defined by an Ideal direction, a lower bound for the step size $t$ can be derived. 
\begin{lemma}\label{second p}
	Suppose that $\mathcal{A}(\bar{{\mathbf x}})\neq\emptyset$. For each $i\in\mathcal{A}(\bar{{\mathbf x}})$ choose $\delta_i$ as in Lemma \ref{zero p} and $\tau, \overline{m}, \overline{W}$ as in Lemma \ref{first p}. Assume that at iteration $k$ of Algorithm \ref{Algorithm 2}, $(\mathbf s_{k1},\ldots,\mathbf s_{k\overline{m}})\in\overline{W}$ and the backtracking Armijo line search in Line 18 is reached with ${\mathbf x}_k\in B(\bar{{\mathbf x}},\min\{\tau,\varepsilon_k/3\})$ and $ {\mathbf d}_k={\mathbf d}^\I_{\varepsilon_k}$. Then $t_k\geq\min\{1,\gamma\varepsilon_k/3\}$.
\end{lemma}
{\it Proof}\,
	For the sake of simplicity in notations, let $m:=\overline{m}$ and  $\varepsilon:=\varepsilon_k$. First, we show that for the sampled points $\{\mathbf s_{k1},\ldots,\mathbf s_{km}\}\cup\{\mathbf s_{k0}\}$, we have $g^\I_{\varepsilon,i}\in G_{\varepsilon,i}(\bar{{\mathbf x}}) $ for all $i\in\{1,\ldots,n\}$. Let $\hat{G}^m_{\varepsilon,i}:=\texttt{co}\{\frac{\partial f}{\partial x_i}(\mathbf s_{ki})\}_{i=1}^m$. Since $(\mathbf s_{k1},\ldots,\mathbf s_{km})\in\overline{W}\subset D^m_{\varepsilon-\bar{\varepsilon}}(\bar{{\mathbf x}})$ with $\bar{\varepsilon}=\tau$, we have $\hat{G}^m_{\varepsilon,i}\subset G_{\varepsilon,i}(\bar{{\mathbf x}})$. Moreover, ${\mathbf x}_k\in B(\bar{{\mathbf x}},\varepsilon/3)$ is a differentiable point and hence $\frac{\partial f}{\partial x_i}({\mathbf x}_k)=\frac{\partial f}{\partial x_i}(\mathbf s_{k0})\in G_{\varepsilon,i}(\bar{{\mathbf x}})$. Thus, $g^\I_{\varepsilon,i}\in G_{\varepsilon,i}(\bar{{\mathbf x}})$ for all $i\in\{1,\ldots,n\}$. Now, set $\delta:=\min\{\delta_i \,:\, i\in\mathcal{A}(\bar{{\mathbf x}})\}$, then the fact that $(\mathbf s_{k1},\ldots,\mathbf s_{km})\in\overline{W}\subset W^m_{\varepsilon}(\bar{{\mathbf x}},\bar{{\mathbf x}},\delta,\omega)$ along with Corollary \ref{corollary zero} imply that
	$$\lvert g^\I_{\varepsilon,i}\rvert < \rho_{\varepsilon,i}(\bar{{\mathbf x}})+\delta  \qquad \text{for \,all} \,\, i\in \mathcal{A}(\bar{{\mathbf x}}).$$ 
	Now by Lemma \ref{zero p}, for each $v_i\in G_{\varepsilon,i}(\bar{{\mathbf x}})$ and $c\in(0,1)$ one can write
	\begin{equation}\label{v1}
	g^\I_{\varepsilon,i} v_i > c\,({g^\I_{\varepsilon,i}})^2 \qquad \text{for \,all} \, i\in \mathcal{A}(\bar{{\mathbf x}}).
	\end{equation}
	Furthermore, Corollary \ref{corollary zero} yields that
	\begin{equation}\label{v2}
	g^\I_{\varepsilon,i}=0, \qquad \text{for\, all}\, i\in\mathcal{A}'_\texttt{int}(\bar{{\mathbf x}}),
	\end{equation}
	and for arbitrary $\omega>0$, we have
	\begin{equation}\label{v3}
	\lvert g^\I_{\varepsilon,i}\rvert \leq \omega,  \qquad \text{for\, all}\, i\in\mathcal{A}'_\texttt{bd}(\bar{{\mathbf x}}).
	\end{equation}
	Using \eqref{v1}, \eqref{v2}, \eqref{v3} and choosing $\omega>0$ sufficiently small in Lemma \ref{first p}, for all $\mathbf v\in G_\varepsilon(\bar{{\mathbf x}})$ we have
	\begin{align*}
	\langle {\mathbf g}^\I_\varepsilon, \mathbf v\rangle &= \sum_{i\in\mathcal{A}(\bar{{\mathbf x}})} g^\I_{\varepsilon,i} v_i + \sum_{i\in\mathcal{A}'_\texttt{int}(\bar{{\mathbf x}})} g^\I_{\varepsilon,i} v_i + \sum_{i\in\mathcal{A}'_\texttt{bd}(\bar{{\mathbf x}})} g^\I_{\varepsilon,i} v_i\\
	&= \sum_{i\in\mathcal{A}(\bar{{\mathbf x}})} g^\I_{\varepsilon,i} v_i + \sum_{i\in\mathcal{A}'_\texttt{bd}(\bar{{\mathbf x}})} g^\I_{\varepsilon,i} v_i > c \sum_{i\in\mathcal{A}(\bar{{\mathbf x}})} ({g^\I_{\varepsilon,i}})^2\\
	&= c \lVert {\mathbf g}^\I_\varepsilon \rVert ^2.
	\end{align*}  
	We stress that, the last inequality follows from \eqref{v1} and the fact that $\omega>0$ is sufficiently small. So far, we have proved that
	\begin{equation}\label{v4}
	\langle {\mathbf g}^\I_\varepsilon, \mathbf v\rangle > c	 \lVert {\mathbf g}^\I_\varepsilon \rVert ^2, \qquad \text{for\, all} \, \mathbf v\in G_\varepsilon(\bar{{\mathbf x}}).	
	\end{equation}
	By contradiction, suppose that $t_k<\min\{1,\gamma\varepsilon/3\}$. In view of Algorithm \ref{Algorithm 1(BLS)}, we conclude
	$$-c\; \gamma^{-1}t_k\lVert {\mathbf g}^\I_\varepsilon\rVert \leq f({\mathbf x}_k+\gamma^{-1}t_k \mathbf d^\I_\varepsilon)-f({\mathbf x}_k), $$
	and using Lebourg's mean value theorem \cite{Clarke1990}, there exist $\hat{{\mathbf x}}_k\in[{\mathbf x}_k+\gamma^{-1}t_k\mathbf d^\I_\varepsilon, {\mathbf x}_k]$ and ${\boldsymbol \xi}_k\in\partial f(\hat{{\mathbf x}}_k)$ such that
	$$f({\mathbf x}_k+\gamma^{-1}t_k\mathbf d^\I_\varepsilon)-f({\mathbf x}_k)=\gamma^{-1}t_k \langle {\boldsymbol \xi}_k,\mathbf d^\I_\varepsilon\rangle. $$
	Therefore, using $\mathbf d^\I_\varepsilon=-\lVert {\mathbf g}^\I_\varepsilon\rVert^{-1}{\mathbf g}^\I_\varepsilon$, we obtain $\langle {\boldsymbol \xi}_k, {\mathbf g}^\I_\varepsilon\rangle \leq c \lVert {\mathbf g}^\I_\varepsilon\rVert^2$ and \eqref{v4} implies that ${\boldsymbol \xi}_k\notin G_{\varepsilon}(\bar{{\mathbf x}})$. On the other hand, $\gamma^{-1}t_k\lVert \mathbf d^\I_\varepsilon\rVert< \varepsilon/3$ and $\lVert {\mathbf x}_k-\bar{{\mathbf x}}\rVert\leq \varepsilon/3$ yield $\hat{{\mathbf x}}_k\in B(\bar{{\mathbf x}},2\varepsilon/3)$ and hence ${\boldsymbol \xi}_k\in G_\varepsilon(\bar{{\mathbf x}})$, which is  a contradiction. \qed

In the next lemma, we recall a highly useful inequality from \cite{Kiwiel2007}.
\begin{lemma}[\cite{Kiwiel2007}]
	Algorithm \ref{Algorithm 2} ensures that 
	\begin{equation}\label{key inequality}
	f({\mathbf x}_{k+1})\leq f({\mathbf x}_k)-\frac{1}{2}c \lVert {\mathbf x}_{k+1}-{\mathbf x}_{k}\rVert \lVert \mathbf g_k\rVert, \qquad {\rm{for\,\, all}} \,\, k,
	\end{equation}
	in which $c$ is the Armijo parameter.
\end{lemma} 

Now, we are in a position to analyze the global convergence behavior of the GSI algorithm. Our analysis is similar to the work of Kiwiel in \cite{Kiwiel2007}. However, there are some differences due to using Ideal directions. So, for the sake of completeness, we present a comprehensive proof for the following theorem.
\begin{theorem}
	Let $\{{\mathbf x}_k\}$ be a sequence generated by Algorithm \ref{Algorithm 2} with $\varepsilon_0, \nu_0>0, \nu_{opt}=\varepsilon_{opt}=0$ and $0<\mu,\theta<1$. Then, with probability 1, the algorithm does not terminate in Line 5 and either $f({\mathbf x}_k)\to-\infty$ or $\nu_k\downarrow0, \varepsilon_{k}\downarrow0$ and every cluster point of the sequence $\{{\mathbf x}_k\}$ is a Clarke stationary point. 
\end{theorem}   
{\it Proof}
	From measure theory, we know that the termination in Line 5 has zero probability. Also, in the case $f({\mathbf x}_k)\to-\infty$, we have nothing to prove.
	So, we consider the following case
	$$\inf_{k\to\infty} f({\mathbf x}_k)>\infty. $$
	Then, sufficient decrease condition \eqref{sufficient decrease} yields
	\begin{equation}\label{THinq1}
	\sum_{k=0}^{\infty}t_k\lVert \mathbf g_k\rVert<\infty
	\end{equation}
	and from \eqref{key inequality}, we have
	\begin{equation}\label{THinq2}
	\sum_{k=0}^{\infty} \lVert {\mathbf x}_{k+1}-{\mathbf x}_k\rVert \lVert \mathbf g_k\rVert < \infty.  
	\end{equation}
	Now, we prove that, with probability 1, the sequences $\{\nu_k\}$ and $\{\varepsilon_k\}$ tend to zero. For this purpose, assume that there exist $k_1\in \mathbb{N},\; \bar{\varepsilon}>0$ and $ \bar{\nu}>0$ such that $\varepsilon_k=\bar{\varepsilon}$ and $\nu_k=\bar{\nu}$ for all $k\geq k_1$. Therefore, we conclude that $\lVert \mathbf g_k\rVert>\bar{\nu}$ in \eqref{THinq1} and \eqref{THinq2} which implies  $t_k\to 0$ and $\lVert {\mathbf x}_{k+1}-{\mathbf x}_k\rVert\to 0$ as $k\to\infty$. Thus, we may assume that there is  $\bar{{\mathbf x}}\in\mathbb{R}^n$ such that ${\mathbf x}_k\to\bar{{\mathbf x}}$. Suppose in the GSI algorithm 
	$$m\geq \overline{m}=2 (\lvert \mathcal{A}(\bar{{\mathbf x}})\rvert+\lvert\mathcal{A}'_\texttt{bd}(\bar{{\mathbf x}})\rvert) +4 \lvert \mathcal{A}'_\texttt{int}(\bar{{\mathbf x}})\rvert \geq n+1. $$
	For the sake of simplicity in notations, let  $\varepsilon:=\bar{\varepsilon}$. We need to consider four cases.\\
	\textit{Case 1.} Let $\mathbf 0\notin G_\varepsilon(\bar{{\mathbf x}})$ and $\mathcal{K}$ be an infinite subset of $\mathbb{N}$ such that $\mathbf d_k=-\lVert {\mathbf g}^\s_{\varepsilon_k}\rVert^{-1} {\mathbf g}^\s_{\varepsilon_k}$ for all $k\geq k_1$ and $k\in\mathcal{K}$. Then, for  $\tau, \delta, $ and $\overline{V}$ chosen as in Lemma \ref{K3}, since $t_k\to 0 $ and ${\mathbf x}_k\to \bar{{\mathbf x}}$ , we can choose $k_2\geq k_1$ such that ${\mathbf x}_k\in B(\bar{{\mathbf x}},\min\{\tau,\varepsilon/3\}) $ and $t_k<\min\{1,\gamma\varepsilon/3\}$ for all $k\geq k_2$ which means that  $(\mathbf s_{k1},\ldots,\mathbf s_{km})\notin \overline{V}$ for all $k\geq k_2$ and $k\in\mathcal{K}$. Since $(\mathbf s_{k1},\ldots,\mathbf s_{km})$ is sampled uniformly and independently from $D^m_\varepsilon({\mathbf x_k})$ containing the nonempty open set $\overline{V}$, this event has probability zero. \\
	\textit{Case 2.} Let $0\notin G_{\varepsilon,i}(\bar{{\mathbf x}})$ for some $i\in\{1,\ldots,n\}$  and $\mathcal{K}$ be an infinite subset of $\mathbb{N}$ such that $\mathbf d_k=-\lVert {\mathbf g}^\I_{\varepsilon_k}\rVert^{-1} {\mathbf g}^\I_{\varepsilon_k}$ for all $k\geq k_1$ and $k\in\mathcal{K}$. Then, for $ \tau, \delta, \omega$ and $\overline{W}$ chosen as in Lemma \ref{second p}, one can choose $k_2\geq k_1$ such that ${\mathbf x}_k\in B(\bar{{\mathbf x}},\min\{\tau,\varepsilon/3\})$ and $t_k<\min\{1,\gamma\varepsilon/3\}$ for all $k\geq k_2$ which implies that $ (\mathbf s_{k1},\ldots,\mathbf s_{km})\notin \overline{W}$ for all $k\geq k_2$ and $k\in\mathcal{K}$. This event occurs with probability zero as well.\\
	\textit{Case 3.} Let $\mathbf 0\in G_{\varepsilon}(\bar{{\mathbf x}})$ and $\mathcal{K}$ be an infinite subset of $\mathbb{N}$ such that $\mathbf d_k=-\lVert {\mathbf g}^\s_{\varepsilon_k}\rVert^{-1} {\mathbf g}^\s_{\varepsilon_k}$ for all $k\geq k_1$ and $k\in\mathcal{K}$. Clearly
	$$\bar{\nu}\leq \lVert {\mathbf g}^\s_{\varepsilon_k} \rVert\leq \text{dist}\Big(\mathbf 0,\texttt{co}\{\nabla f(\mathbf s_{k1}),\ldots,\nabla f(\mathbf s_{km})\}\Big). $$
	Then, for $\delta:=\bar{\nu}/2$ and $\tau, \overline{V}$ as in Lemma \ref{K2}, one can choose  $k_2\geq k_1$ such that ${\mathbf x}_k\in B(\bar{{\mathbf x}},\tau)$ for all $k\geq k_2$ and $k\in\mathcal{K}$ and since $\rho_\varepsilon(\bar{{\mathbf x}})=0$, we can write
	$$\text{dist}\Big(\mathbf 0,\texttt{co}\{\nabla f(\mathbf s_{k1}),\ldots,\nabla f(\mathbf s_{km})\}\Big)>\rho_\varepsilon(\bar{{\mathbf x}})+\bar{\nu}/2=\bar{\nu}/2, $$
	which means that $(\mathbf s_{k1},\ldots,\mathbf s_{km})\notin \overline{V}$ for all $k\geq k_2$ and $k\in\mathcal{K}$. This event takes place with probability zero.\\
	\textit{Case 4.} Let $0\in G_{\varepsilon,i}(\bar{{\mathbf x}})$ for every $i\in\{1,\ldots,n\}$ and $\mathcal{K}$ be an infinite subset of $\mathbb{N}$ such that $\mathbf d_k=-\lVert {\mathbf g}^\I_{\varepsilon_k}\rVert^{-1} {\mathbf g}^\I_{\varepsilon_k}$ for all $k\geq k_1$ and $k\in\mathcal{K}$. If $\mathcal{A}'_\texttt{bd}(\bar{{\mathbf x}})=\emptyset$, then for $\tau, \overline{W}$ as in Lemma \ref{first p}, one can choose  $k_2\geq k_1$ such that ${\mathbf x}_k\in B(\bar{{\mathbf x}},\tau)$ for all $k\geq k_2$ and the fact that
	$$0<\bar{\nu}\leq\lVert {\mathbf g}^\I_{\varepsilon_k} \rVert, \qquad \text{for \,all} \,\, k\geq k_2 \,\, \text{and} \,\, k\in\mathcal{K}, $$
	gives $(\mathbf s_{k1},\ldots,\mathbf s_{km})\notin\overline{W}$ for all $k\geq k_2$ and $k\in\mathcal{K}$. Clearly, this event occurs with probability zero. Now, suppose that $\mathcal{A}'_\texttt{bd}(\bar{{\mathbf x}})\neq\emptyset$. Then, for $\omega:=\lvert \mathcal{A}'_\texttt{bd}(\bar{{\mathbf x}})\rvert^{-1/2}\bar{\nu}/2$
	and $\tau,\overline{W}$ as in Lemma \ref{first p}, one can choose  $k_2\geq k_1$ such that ${\mathbf x}_k\in B(\bar{{\mathbf x}},\tau)$ for all $k\geq k_2$ and the fact that
	$$0<\omega<\bar{\nu}\leq\lVert {\mathbf g}^\I_{\varepsilon_k}\rVert, \qquad \text{for \,all} \,\, k\geq k_2 \,\, \text{and} \,\, k\in\mathcal{K},  $$
	yields $(\mathbf s_{k1},\ldots,\mathbf s_{km})\notin\overline{W}$ for all $k\geq k_2$ and $k\in\mathcal{K}$. This event has also probability zero. \\
	Now, we are in a position to consider the case that $\nu_k\downarrow0, \varepsilon_k\downarrow0$ and $\bar{{\mathbf x}}$ is a cluster point of $\{{\mathbf x}_k\}$. Since $\nu_k\downarrow0$ we conclude that $0$ is a cluster point of the sequence $\{\lVert {\mathbf g}^\s_{\varepsilon_k}\rVert\}$. Therefore, in the case ${\mathbf x}_k\to\bar{{\mathbf x}}$ we have
	$$\liminf_{k\to\infty} \max\{\lVert {\mathbf x}_k-\bar{{\mathbf x}}\rVert,\lVert {\mathbf g}^\s_{\varepsilon_k}\rVert, \varepsilon_k\}=0. $$
	Note that, ${\mathbf g}^\s_{\varepsilon_k}\in\partial_{\varepsilon_k}f({\mathbf x}_k)$ and hence outer semicontinuity of the map $\partial_\cdot f(\cdot)$ yields $\mathbf 0\in\partial f(\bar{{\mathbf x}}).$ On the other hand, if ${\mathbf x}_k\nrightarrow\bar{{\mathbf x}}$ we need to prove that
	\begin{equation}\label{THinq4}
	\liminf_{k\to\infty} \max\{\lVert {\mathbf x}_k-\bar{{\mathbf x}}\rVert,\lVert {\mathbf g}^\s_{\varepsilon_k}\rVert\}=0. 
	\end{equation}
	By contradiction, suppose that there exist $\bar{\nu}>0, k_1\in\mathbb{N}$ and an infinite set $\mathcal{K}:=\{k \,:\, k\geq k_1,  \lVert {\mathbf x}_k-\bar{{\mathbf x}}\rVert\leq\bar{\nu}, \lVert {\mathbf g}^\s_{\varepsilon_k}\rVert>\bar{\nu}\}$. Then, by \eqref{THinq2} we have
	\begin{equation}\label{THinq3}
	\sum_{k\in\mathcal{K}}\lVert {\mathbf x}_{k+1}-{\mathbf x}_k\rVert<\infty .
	\end{equation}
	By noting that ${\mathbf x}_k\nrightarrow\bar{{\mathbf x}}$, there is $\varepsilon>0$ such that for all $k\in\mathcal{K}$ with $\lVert {\mathbf x}_k-\bar{{\mathbf x}}\rVert\leq\bar{\nu}/2$ there exists $k_2>k$ satisfying $\lVert {\mathbf x}_{k_2}-{\mathbf x}_k\rVert>\varepsilon$ and $\lVert {\mathbf x}_i-\bar{{\mathbf x}}\rVert\leq\bar{\nu}$ for each $k\leq i\leq k_2$. Therefore, using the triangle inequality, we have $\varepsilon<\lVert {\mathbf x}_{k_2}-{\mathbf x}_k\rVert\leq \sum_{i=k}^{k_2-1}\lVert {\mathbf x}_{i+1}-{\mathbf x}_i\rVert$. However, for $k\in\mathcal{K}$ sufficiently large, \eqref{THinq3} yields that the right-hand side of
	this inequality is less than $\varepsilon$, which is a contradiction. Therefore \eqref{THinq4} holds and consequently $\mathbf 0\in\partial f(\bar{{\mathbf x}})$. \qed

\subsection{More results }
In this subsection, we modify the Algorithm \ref{Algorithm 2}, such that under some moderate conditions on the the objective function $f$, one can obtain an upper bound for the number of serious iterations required for obtaining a $(\nu, G^m_\varepsilon({\mathbf x}_k))$-stationary point. Next, the convergence of the proposed method is studied. Indeed, instead of using the backtracking Armijo line search presented in Algorithm \ref{Algorithm 1(BLS)}, following \cite{Kiwiel2007}, we will apply the following limited backtracking Armijo line search.  

\LinesNumbered
\IncMargin{1em}
\begin{algorithm}[H]
	\SetAlgoLined
	
	\SetKwInOut{Input}{inputs}
	\SetKwInOut{Output}{output}
	\SetKwInOut{Required}{requied}
	
	\Indm  
	\Input{${\mathbf x}\in\mathbb{R}^n$, the directions ${\mathbf d},\mathbf g\in \mathbb R^n$, $\varepsilon>0$, parameters   $\gamma, c\in(0,1)$. }
	\Output{Step size $t$.}
	\Required{${\mathbf d}=-\lVert \mathbf{g}\rVert^{-1} \mathbf g$ is an ADD for $f$ at ${\mathbf x}$.}
	\Indp
	\BlankLine
	
	\SetKwFunction{FMain}{LBALS}
	\SetKwProg{Fn}{Function}{:}{}
	\Fn{\FMain{$\mathbf x$, $\mathbf g$, $\mathbf d$, $\varepsilon$, $\gamma$, $c$}}{
		$t:=1$\;
		\While{$t>\min\{1,\gamma\varepsilon/3\}$}{
			\If{$f({\mathbf x}+t{\mathbf d})-f({\mathbf x})<-c t \lVert {\mathbf g} \rVert$}{\textbf{return} $t$\;} 
			$t:=\gamma t$\;
			
		}
		\textbf{return} $t:=0$\; 
	}
	\textbf{End Function}
	\caption{Limited Backtracking Armijo Line Search (LBALS) \cite{Kiwiel2007} }
	\label{Algorithm 3(LBALS)}
\end{algorithm}
\DecMargin{1em}

By this substitution, in some iterations  of  the Algorithm \ref{Algorithm 2}, the returned step size $t$ by LBALS may be zero. In theses iterations, which are called ``null iterations'',  no improvement is achieved. When a null iteration takes place, then in the next iteration(s),  the gradient bundle $G^m_{\varepsilon_k}({\mathbf x}_k)$ is resampled, such that eventually a serious iteration ($t_k\ne 0$ )  occurs when  the resulting gradient bundle improves the search direction sufficiently.

At the first glance, the null iterations are undesirable, since we do some computations without any reduction in the objective function. However, according to \cite{Kiwiel2007}, this change has some advantages in reducing the number of function evaluations. Moreover, as we will show, it does not affect the convergence of the  GSI method and 
interesting results can be explored. Maybe the most important of them is the subject of the following theorem in which the LBALS is considered as an alternative of BALS.

\begin{theorem}\label{TH-1}
	In the GSI algorithm, let $ \nu_0=\nu_{opt}=:\nu>0, \varepsilon_0=\varepsilon_{opt}=:\varepsilon>0$ and $m:=4n$. If the function $f$ is bounded below, then the algorithm terminates after finitely iterations and produces a $(\nu,G^m_\varepsilon({\mathbf x}_k))$-stationary point, with probability 1.
\end{theorem}
{\it Proof}
	We know that the GSI algorithm does not terminate in Line 5, with probability 1. So, we consider the case that the Algorithm \ref{Algorithm 2} does not stop in Line 5. By contradiction, suppose that an infinite sequence $\{{\mathbf x}_k\}$ is generated by this algorithm. Then, with probability 1, there exists an infinite set $\mathcal{K}$ such that $t_k\neq0$ for all $k\in\mathcal{K}$. Let $\mathcal{K}=\mathcal{K}_1\cup\mathcal{K}_2$ such that $\lVert {\mathbf g}^\I_{\varepsilon_k}\rVert\leq\nu$ at any iteration $k\in\mathcal{K}_1$ and $\lVert {\mathbf g}^\I_{\varepsilon_k}\rVert>\nu$ at any iteration $k\in\mathcal{K}_2$. For every $k\in\mathcal{K}_2$, the direction $\mathbf d_k=-\lVert {\mathbf g}^\I_{\varepsilon_k}\rVert^{-1} {\mathbf g}^\I_{\varepsilon_k}$ is considered as the search direction and we have
	$$f({\mathbf x}_k+t_k\mathbf d_k)-f({\mathbf x}_k)\leq -ct_k\lVert{\mathbf g}^\I_{\varepsilon_k}\rVert, \qquad \text{for \,all} \,\, k\in\mathcal{K}_2. $$
	In view of condition $\lVert {\mathbf g}^\I_{\varepsilon_k}\rVert>\nu$, one can write
	\begin{equation}\label{THinq5}
	f({\mathbf x}_{k+1})-f({\mathbf x}_k)\leq -ct_k\nu,  \qquad \text{for \,all} \, k\in\mathcal{K}_2.  
	\end{equation}
	On the other hand, for every $k\in\mathcal{K}_1$, since the algorithm does not terminate, we have $\lVert {\mathbf g}^\s_{\varepsilon_k}\rVert>\nu$ for all $k\in\mathcal{K}_1$ and similarly, we can write
	\begin{equation}\label{THinq6}
	f({\mathbf x}_{k+1})-f({\mathbf x}_k)\leq -ct_k\nu,\quad \text{for \,all} \,\, k\in\mathcal{K}_1.  
	\end{equation}
	Combining \eqref{THinq5} with \eqref{THinq6} we obtain
	$$ f({\mathbf x}_{k+1})-f({\mathbf x}_k)\leq -ct_k \nu , \quad \text{for \,all} \,\, k\in \mathcal{K} .$$
	Furthermore, for every $k\in\mathbb{N}_0\setminus\mathcal{K}$, we have $t_k=0$ and hence
	$$f({\mathbf x}_{k+1})-f({\mathbf x}_k)=-c\nu t_k=0, \quad \text{for \,all} \,\, k\in\mathbb{N}_0\setminus\mathcal{K}, $$
	and consequently
	$$ f({\mathbf x}_{k+1})-f({\mathbf x}_k)\leq -ct_k \nu , \quad \text{for \,all} \,\, k.$$
	Using the above inequality inductively, one can write
	\begin{equation}\label{THinq7}
	f({\mathbf x}_{k+1})\leq f({\mathbf x}_0)-c\nu \sum_{j=0}^{k}t_j ,\quad \text{for \,all} \,\, k.
	\end{equation}
	Since $t_j> \min\{1,\gamma\varepsilon/3\}$ for every $j\in\mathcal{K}$, we have 
	$$\sum_{j\in\mathcal{K} }t_j=\infty,$$
	and hence the left hand side of the inequality \eqref{THinq7} tends to $-\infty$ as $k\to\infty$, while function $f$ is supposed to be bounded from below.\qed

\begin{corollary}
	Under assumptions of Theorem \ref{TH-1}, with probability 1, the GSI Algorithm terminates after finitely serious iterations $k_{\max}\geq0$ such that
	$$ k_{\max}\leq \Big\lfloor \, \frac{f(\mathbf x_0)-f_l}{c\nu \min\{1,\gamma\varepsilon/3\}} \,\Big\rfloor +1 ,$$
	where $f_l$ is a lower bound for $f$. 	
\end{corollary}
{\it Proof}
	By contradiction, suppose that an infinite sequence $\{{\mathbf x}_k\}$ is generated by this algorithm. Similar to the proof of Theorem \ref{TH-1}, we have 
	$$f({\mathbf x}_{k+1})-f({\mathbf x}_k)\leq-c t_k\nu ,\quad \text{for \,all} \,\, k,$$
	and therefore
	$$f({\mathbf x}_{k+1})\leq f({\mathbf x}_0)-K c\, \nu \min\{1,\gamma\varepsilon/3\}, \quad \text{for \,all} \,\, k, $$
	in which
	$$K:=\Big\lvert \{j\in\mathcal{K} \, : \, j\leq k \}\Big\rvert. $$
	Since $K\to\infty$ as $k\to\infty$, we conclude that $f({\mathbf x}_k)$ tends to $-\infty$ as $k\to\infty$ which is a contradiction.
	Moreover, since $f_l$ is a lower bound for the function $f$, one can write
	$$ f_l\leq f({\mathbf x}_{k+1})\leq f({\mathbf x}_0)-Kc\nu \min\{1,\gamma\varepsilon/3\} .$$
	Hence, an upper bound  for the number of serious iterations is obtained as follows
	$$ k_{\max}\leq \Big\lfloor \frac{f({\mathbf x}_0)-f_l}{c\nu \min\{1,\gamma\varepsilon/3\}}\Big\rfloor +1 .$$\qed
	
In the next theorem, we examine the convergence analysis of Algorithm \ref{Algorithm 2} under the proposed modifications. In this theorem,
the lower $\alpha$-level set of a function $f:\mathbb{R}^n\to\mathbb{R}$ is given by \cite{Rockafellar2004}
$$\text{lev}_{\leq\alpha}f:=\{{\mathbf x}\in\mathbb{R}^n \,\,: \,\, f({\mathbf x})\leq\alpha\}  .$$

\begin{theorem}\label{TH-2}
	Suppose that the objective function $f$ is bounded from below and  $\text{lev}_{\leq f({\mathbf x}_0)}f$ is bounded. In Algorithm \ref{Algorithm 2}, let $\varepsilon_0, \nu_0>0, \nu_{opt}=\varepsilon_{opt}=0$ and $0<\mu,\theta<1$. Then, with probability 1, any cluster point of the sequence $\{{\mathbf x}_k\}$ generated by this algorithm is a Clarke stationary point.
\end{theorem} 
{\it Proof}
	According to Theorem \ref{TH-1}, with probability 1, the GSI Algorithm produces a sequence of $(\nu_i, G^m_{\varepsilon_i}({\mathbf x}_i))$-stationary points, say $\{{\mathbf x}_i\}$. Hence
	$$ \min \{\lVert \mathbf g \rVert \,\, : \,\, {\mathbf g}\in  G^m_{\varepsilon_i}({\mathbf x}_i) \}\leq\nu_i.$$
	Clearly, ${\mathbf x}_i\in \text{lev}_{\leq f(\mathbf x_0)}f$ for every $i$. Thus, the boundedness of  $\text{lev}_{\leq f({\mathbf x}_0)}f$ implies that the sequence $\{{\mathbf x}_i\}$ has at least one cluster point. So, assume that $\{{\mathbf x}_{i_j}\}\subset \{{\mathbf x}_i\}$ be such that ${\mathbf x}_{i_j}\to {\bar{\mathbf x}}$ as $j\to\infty$. Therefore
	\begin{align}\label{A1}
	\min \{\lVert \mathbf g \rVert \,\, : \,\, {\mathbf g}\in  G^m_{\varepsilon_{i_j}}({\mathbf x}_{i_j}) \}\leq\nu_{i_j}.
	\end{align}
	Now, let $\omega>0$ be arbitrary. Since $\nu_i\to 0$, there exists $j_1>0$ such that $\nu_{i_j}<\omega$ for all $j>j_1$. In view of \eqref{A1} we have
	$$\min \{\lVert {\mathbf g} \rVert \,\, : \,\, {\mathbf g}\in  G^m_{\varepsilon_{i_j}}({\mathbf x}_{i_j}) \}<\omega, \quad \text{for \,all} \,\, j>j_1 .$$
	Now, we have
	$$\lVert {\mathbf g}^\s_{\varepsilon_{i_j}}\rVert =\min \{\lVert {\mathbf g} \rVert \,\, : \,\, {\mathbf g}\in  G^m_{\varepsilon_{i_j}}({\mathbf x}_{i_j}) \},  \quad \text{for \,all} \,\, j>j_1,$$
	and since $\lVert {\mathbf g}^\s_{\varepsilon_{i_j}}\rVert<\omega$ for every $j>j_1$, we may assume that $ {\mathbf g}^\s_{\varepsilon_{i_j}}\to {\bar{\mathbf g}}$ as $j\to\infty$. On the other hand, for every $j>j_1$ we have ${\mathbf g}^\s_{\varepsilon_{i_j}}\in  G^m_{\varepsilon_{i_j}}(\mathbf x_{i_j})$ and therefore $ {\mathbf g}^\s_{\varepsilon_{i_j}}\in \partial_{\varepsilon_{i_j}}f({\mathbf x}_{i_j})$. Now, using outer semicontinuity of the set valued map $\partial_\cdot f(\cdot)$, we conclude
	$$ {\bar{\mathbf g}}\in \partial f({\bar{\mathbf x}}) ,$$
	which yields
	$$  \min \{\lVert {\boldsymbol \xi} \rVert \,\, : \,\, {\boldsymbol \xi} \in  \partial f({\bar{\mathbf x}}) \}<\omega .$$
	Since $\omega>0$ is arbitrary we conclude that
	$$ \mathbf 0\in\partial f({\bar{\mathbf x}}).$$  \qed
	According to the results of this subsection, when the function $f$ is bounded below, using  LBALS instead of BALS guarantees an upper bound for the number of serious iterations required for obtaining a $(\nu, G^m_\varepsilon(\mathbf x_k))$-stationary point.

\section{Numerical results}\label{Section6}
In this section, we demonstrate the efficiency of the proposed method through the results of numerical experiments. Following \cite{Bagirov2014},  we divide test problems into three categories based on the scale of the problems in which problems with $n\leq50$ are considered as small scale problems, medium scale problems have a dimension $50< n\leq 200$ and large scale problems are identified by $n>200$.

 Based on this classification, this section is divided into three experiments. In Experiment 1, we apply the GSI method to a set of nonsmooth small scale problems. Medium and large scale problems are considered in Experiments 2 and 3, where performance profiles are provided to compare the proposed method with the GS method. At each experiment, a set of nonsmooth test problems including convex and nonconvex objectives is considered.

The experiments are implemented in \textsc{Matlab} software on PC Intel Core i7 2700k CPU 3.5 GHz and 8 GB of RAM. In Experiment 3, we also used \textsc{Matlab Api} for C++ language in order to manage the CPU time elapsed by each run. To provide fair results,  we run each problem 5 times using starting points generated randomly from a ball centered at ${\mathbf x}_0$ (suggested in the literature) with radius $\lVert {\mathbf x}_0\lVert/n$. Note that, ${\mathbf x}_0\neq \mathbf 0$ for all problems, so the starting points for each run were different.  In the test problems for which ${\mathbf x}_0$ is not available in the literature, we used ${\mathbf x}_0=\mathbf e:=(1,1,\ldots,1)$ when ${\mathbf x}^*\neq \mathbf e$ and ${\mathbf x}_0=2 \mathbf e$ when ${\mathbf x}^*=\mathbf e$.

For the GS method, we used the default setting proposed in \cite{Burke2005} and for the GSI method we have the following choices for the parameters.

\begin{itemize}
	\item \textit{Sample size}. We have experimented different sizes of sample points. Although Lemma \ref{first p} suggests a sample size $m\geq\bar{m}$, we don't see any favorable change in search direction by choosing $m>2n$.
	 Following \cite{Burke2005}, for the sample size $m=2n$, a suitable approximation of  $\partial_\varepsilon f({\mathbf x})$ is obtained. Since both GS and GSI methods need a great deal of random access memory, it is not recommended to choose a sample size more than $2n$, specially for large scale objectives. 
	\item \textit{Sampling radius}. We choose $\varepsilon_0=10^{-3}$ for problems with $n\leq 10$ and  $\varepsilon_0=10^{-2}$ for the problems having $n>10$. Furthermore, we set $\mu=0.5$. Note that, there is a trade-off between sampling radius and the efficiency of the Ideal direction. Although a larger sampling radius leads to collecting more information of nonsmooth curves, it may cause ${\mathbf g}^\I_\varepsilon$ to become zero. The GSI has its best performance when there is a balance between these two factors. 
	\item \textit{Backtracking Armijo line search}. In Algorithm \ref{Algorithm 1(BLS)}, we set $\gamma=0.5$, $c=10^{-6}$ and we start with $t=1$. Following \cite{Burke2005}, We also limit the number of backtracking steps to $50$. Once this limit is reached, the sampling radius is reduced by reduction factor $\mu$ and we set ${\mathbf x}_{k+1}:={\mathbf x}_k$, and then we go to Line 3. So, a null iteration occurs and the gradient bundle is resampled with a smaller sampling radius.
	
	\item \textit{Optimality tolerance}. We select $\nu_0$  with respect to the scale of the problems. For small, medium and large scale problems, we set $\nu_0=10^{-3}$, $\nu_0=10^{-2}$ and $\nu_0=10^{-1}$ respectively. Furthermore, we set $\theta=0.5$.
	\item \textit{Stopping criterion}. Since the (local) minimizers of the objective functions are available in our test problems, we terminate the GS and GSI algorithms, if they find a solution which fulfills the following condition   
	
	$$\frac{\lvert f({\mathbf x}_k)-f({\mathbf x}^*)\rvert}{\lvert f({\mathbf x}^*) \rvert+1} < \epsilon. $$
	Here, $\epsilon>0$ is a given tolerance and $f({\mathbf x}^*)$ is the known (local) minimum of the objective function. We set $\epsilon=5\times10^{-4}$ for small scale problems and $\epsilon=10^{-3}$ for  medium and large scale problems. We also limit the number of iterations to $2000$.

\end{itemize}

\subsection{Experiment 1. Small scale problems}\label{exp1}

In this subsection, we apply the proposed method to a set of small scale problems. A description of the considered set of test problems is presented in Table \ref{Table1}.

\begin{table}[ht]
	\centering
	 \caption{List of test problems for Experiment 1.}\label{Table1}

	\begin{tabular}{ c c c c c  } 
		
		\toprule
			Problem&Name & Function Type & $n$ & Ref.  \\
		\hline \\
		1& QL & Convex & $2$ & \cite{Makela_book} \\
		2& Wong1 & Convex & $7$ &  \cite{techreport} \\ 
		3& Wolfe & Convex & $2$ &  \cite{techreport} \\ 
		4& SPIRAL & Convex & $2$&  \cite{techreport}  \\ 
		5& Rosenbrock Function & Nonconvex& $2$ & \cite{Skaaja} \\
		6& Crescent & Nonconvex & $2$ & \cite{kiwielbook}  \\
		7& Mifflin2& Nonconvex & $2$ &  \cite{Makela_book} \\
		8& EVD52 & Nonconvex & $3$  & \cite{techreport} \\
		9& HS78  & Nonconvex & $5$ & \cite{techreport}  \\
		10& Condition Number & Nonconvex & $45$ & \cite{Skaaja}  \\
		
		\bottomrule
		
	\end{tabular}
	
\end{table}

\begin{table}[ht]
	
	\centering
	\caption{Numerical results for  small scale problems.} \label{Table2}
	\begin{tabular}{cccccc} %
		\toprule
		& \multicolumn{5}{c}{GSI}\\
		\cmidrule{2-6}
		Problem & Iters & NII & PII & $f_{\text{eval}}$ & $g_{\text{eval}}$\\
		\cmidrule{1-6}
		
		1 & 14  & 12 & 86 \%   & 76   & 70    \\
		2 & 82  & 75 & 91 \%   & 730  & 1233  \\
		3 & 29  & 29 & 100\%   & 172  & 145   \\
		4 & 1735&1712& 99 \%   &16390 & 8676  \\
		5 &114  & 76 & 67 \%   & 918  & 568   \\
		6 & 22  & 18 & 82 \%   & 130  & 109   \\
		7 & 7   & 7  & 100\%   & 40   & 36    \\
		8 & 28  & 25 & 89 \%   & 207  & 199   \\
		9 & 85  & 51 & 60 \%   & 1043 & 937   \\
		10& 20  & 20 & 100\%  & 43  & 1801   \\
		
		\bottomrule
	\end{tabular}

\end{table}

Numerical results for this set of test  problems are given in Table \ref{Table2}. In this table, by the term ``Iters'' we refer to the number of iterations. Furthermore, the number of iterations, in which the search direction is obtained by Ideal directions and not by the QP \eqref{CQP}, is denoted by ``NII'' (Number of Ideal Iterations) and the corresponding percentage is abbreviated by ``PII'' (Percentage of Ideal Iterations). In fact, in such iterations the GSI method  only needs to compute the Ideal direction. Also, the number of function and gradient evaluations are denoted by ``$f_{\text{eval}}$'' and ``$g_{\text{eval}}$'' respectively. Since the CPU time elapsed by the GSI method is almost zero for most of these problems, we don't report this factor in this experiment.

It can be seen from Table \ref{Table2} that, the proposed method has been able to solve all the test problems without observing any failure, which means that the GSI method is reliable for this type of problems. In addition, the method has a great potential to find a descent direction without solving the QP \eqref{CQP}. In particular, we see that, for Problems 3, 7 and 10 the method reached its termination tolerance without solving any QP.

\subsection{Experiment 2. Medium scale problems}\label{exp2}
In this experiment, we apply the GS and GSI methods to a set of medium scale problems given in Table \ref{Table3}. All of these problems can be formulated with any number of variables. 

\begin{table}[ht]
	\centering	
	\caption{List of test problems for Experiment 2.}\label{Table3}
	\begin{tabular}{ c c c c   } 
			\toprule 
			Problem&Name & Function Type &  Ref.  \\
		\hline \\

		1& Generalization of L1HILB &Convex & \cite{techreport2}   \\
		2& Generalization of MXHILB & Convex & \cite{Haarala}  \\
		3& Chained LQ & Convex & \cite{Haarala}  \\
		4& Chained CB3 I & Convex & \cite{Haarala} \\
		5& Chained CB3 II & Convex & \cite{Haarala}  \\
		6& Number of Active Faces & Nonconvex & \cite{phdthesis}  \\
		7& Generalization of Brown Function 2 & Nonconvex & \cite{Haarala}\\
		8& Chained Mifflin 2 & Nonconvex & \cite{Haarala}  \\  
		9& Chained Crescent I & Nonconvex &  \cite{Haarala} \\
		10& Chained Crescent II & Nonconvex &  \cite{Haarala} \\
		
		\bottomrule
		\end{tabular}
	
\end{table}

\begin{figure}[H]
	\centering
	\includegraphics[width=\textwidth]{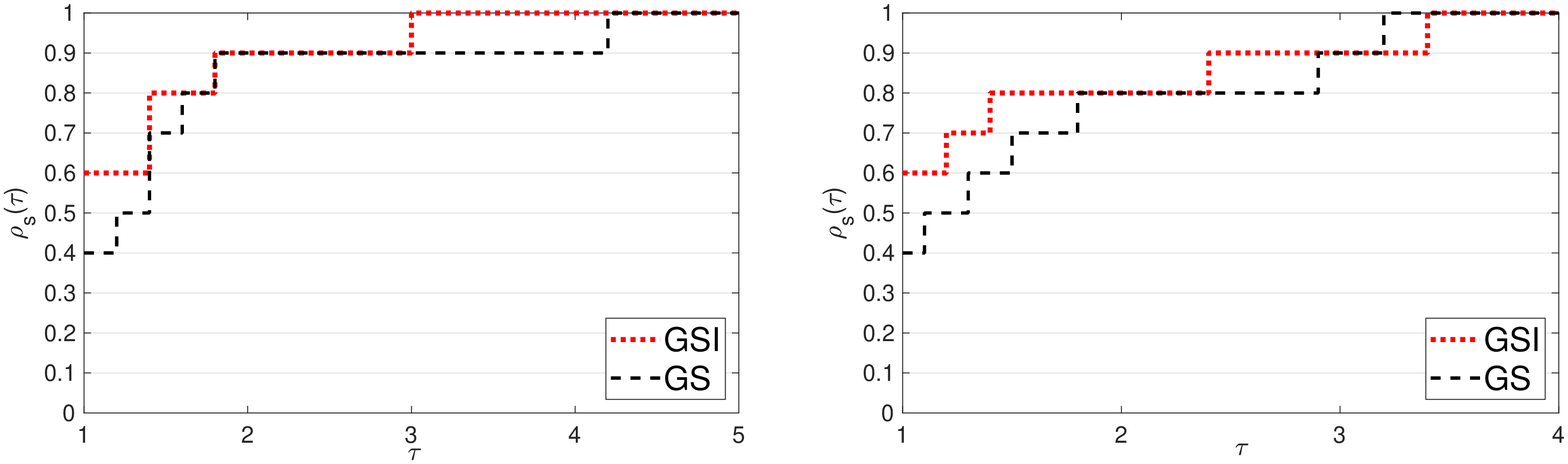}
	\caption{Performance profiles based on CPU time for medium scale problems with $n=100$ (left) and $n=200$ (right). }
	\label{Fig1}
\end{figure}

\begin{figure}[H]
	\centering
	\includegraphics[width=\textwidth]{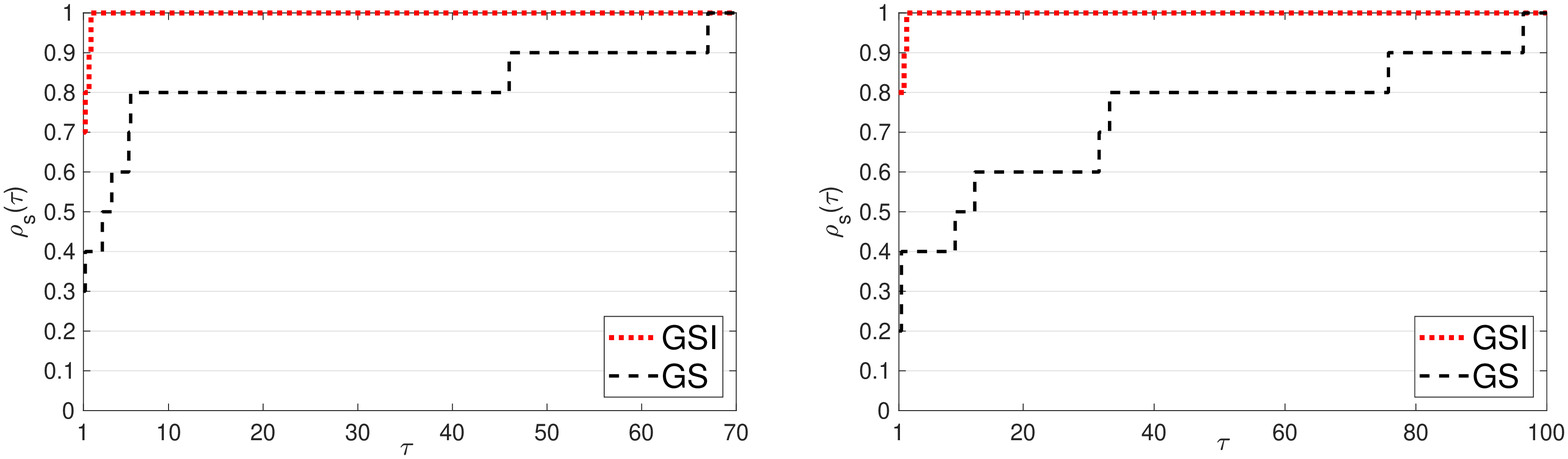}
	\caption{Performance profiles based on QP solver time for medium scale problems with $n=100$ (left) and $n=200$ (right). }
	\label{Fig2}
\end{figure}

We compare our method with the GS method by means of performance profiles \cite{Dolan-More}. 
  Fig. \ref{Fig1} and  Fig. \ref{Fig2} demonstrate results for problems with $n=100$ and $n=200$. It can be seen from Fig. \ref{Fig1} that, for $60\%$ of the problems the proposed method has been more successful than the GS method in the sense of using CPU time. Note that, the efficiency of the GSI method depends on the number of Ideal iterations. Therefore, the GSI method is more successful than the GS method as long as Ideal directions define the search direction in majority of the iterations. Moreover, the efficiency of the Ideal iterations in reducing the CPU time elapsed by the QP solver can obviously be seen in Fig \ref{Fig2}.

\subsection{Experiment 3. Large scale problems}\label{exp3}

In this experiment, we consider a set of nonsmooth test problems for which we set $n=500$ and $n=1000$. A description of the test problems is presented in Table \ref{Table5}.  To the best of our knowledge, it is the first time that the GS method is applied to a set of  problems having more than $200$ variables. Note that, due to heavy computational cost, an unreasonable amount of CPU time may be required for each run. Accordingly, we used \textsc{Matlab Api} for C++ language to manage the CPU time  for the considered large scale problems. 

\begin{figure}[H]
	\centering
	\includegraphics[width=\textwidth]{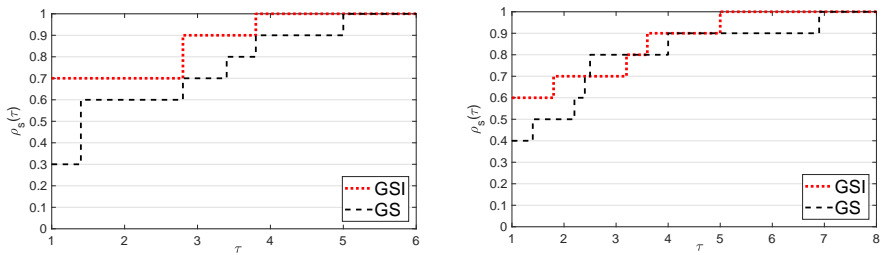}
	\caption{Performance profiles based on CPU time for large scale problems with $n=500$ (left) and $n=1000$ (right). }
	\label{Fig3}
\end{figure}

\begin{figure}[H]
	\centering
	\includegraphics[width=\textwidth]{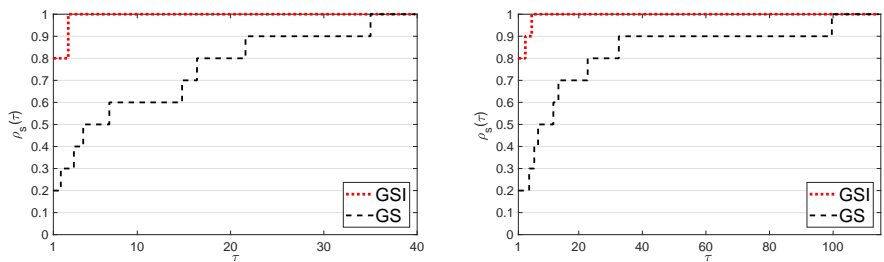}
	\caption{Performance profiles based on QP solver time for large scale problems with $n=500$ (left) and $n=1000$ (right). }
	\label{Fig4}
\end{figure}

Fig. \ref{Fig3} and Fig. \ref{Fig4} show results for problems with  $500$ and $1000$ variables. As Fig. \ref{Fig3} demonstrates, at $70\%$ of the problems having $500$ variables, the GSI method used less CPU time than the GS method. In addition, Fig. \ref{Fig3} also demonstrates  the superiority of the GSI method over the GS method for problems with $1000$ variables. It can be seen from Fig. \ref{Fig4} that, the proposed method has a remarkable potential to reduce the CPU time consumed by the QP solver. Thus, in the problems for which the QP \eqref{CQP} is a demanding part of the GS method, the GSI method is a more desirable choice.

\begin{table}[ht]
	\centering	
	\caption{List of test problems for Experiment 3.}\label{Table5}
	\begin{tabular}{ c c c c   }

		\toprule	Problem&Name & Function Type &  Ref.  \\
		\hline \\

		1& Titled Norm Function &Convex & \cite{BFGS}   \\
		2& A Convex Partly Smooth Function & Convex & \cite{BFGS}  \\
		3& MAXQ & Convex & \cite{Haarala}  \\
		4& Chained LQ & Convex & \cite{Haarala} \\
		5& Nesterov's Chebyshev-Rosenbrock Function1 & Nonconvex & \cite{BFGS}  \\
		6& Nesterov's Chebyshev-Rosenbrock Function2 & Nonconvex & \cite{BFGS}  \\
		
		7& Chained Mifflin 2 & Nonconvex & \cite{Haarala}  \\  
		8& Chained Crescent I & Nonconvex &  \cite{Haarala} \\
		9& Chained Crescent II & Nonconvex &  \cite{Haarala} \\
		10&A Nonconvex Partly Smooth Function.  & Nonconvex &  \cite{BFGS} \\	
		\bottomrule
	\end{tabular}
	
\end{table}

\section{Conclusion}\label{Conclusion}

In this article, by introducing the Ideal direction as an alternative of the approximate $\epsilon$-steepest descent direction, we reduced the need of solving quadratic subproblem in the GS method.  We have shown that, Ideal directions are easy to compute and a nonzero Ideal direction can make a substantial reduction in the objective function. Furthermore, we showed that by using Ideal directions, not only is there no need to consider the quadratic subproblem in the smooth regions, but also it can be an efficient search direction when the method starts to track a nonsmooth curve towards a stationary point. We studied the convergence of the proposed method and we proved that using Ideal directions preserves the global convergence property of the GS method. Furthermore, by using limited backtracking Armijo line search and under some moderate assumptions, we proposed an upper bound for the number of serious iterations in the GSI algorithm. We presented numerical results using small, medium and large scale nonsmooth test problems. Our numerical results clearly demonstrate that, the GSI method inherits robustness from the GS method and it is a significant enhancement in reducing the number of quadratic subproblems.
\begin{acknowledgements}
	The authors would like to thank Dr. L. E. A. Sim\~{o}es and Dr. Frank. E. Curtis for their valuable comments that improved the quality of the paper. 
\end{acknowledgements}


\end{document}